\documentclass[12pt]{article}
\usepackage[]{amsmath,amssymb}
\usepackage{amscd}
\usepackage{latexsym}
\usepackage{cite}

\newtheorem{definition}{Definition}[section]
\newtheorem{theorem}[definition]{Theorem}
\newtheorem{lemma}[definition]{Lemma}
\newtheorem{corollary}[definition]{Corollary}
\newtheorem{remark}[definition]{Remark}

\newtheorem{proposition}[definition]{Proposition}
\newtheorem{notation}[definition]{Notation}

\def\I{\mathbb I}

\def\C{\mathbb C}

\def\K{\mathbb K}

\begin{document}

\title{Hypercubes, Leonard triples and \\ the anticommutator spin algebra}
\author{George M. F. Brown}
\maketitle

\begin{abstract}
This paper is about three classes of objects: Leonard triples, distance-regular graphs and the finite-dimensional irreducible modules for the anticommutator spin algebra.  Let $\K$ denote an algebraically closed field of characteristic zero.  Let $V$ denote a vector space over $\K$ with finite positive dimension.  A Leonard triple on $V$ is an ordered triple of linear transformations in $\mathrm{End}(V)$ such that for each of these transformations there exists a basis for $V$ with respect to which the matrix representing that transformation is diagonal and the matrices representing the other two transformations are irreducible tridiagonal.  The Leonard triples of interest to us are said to be totally B/AB and of Bannai/Ito type.

It is known that totally B/AB Leonard triples of Bannai/Ito type arise in conjunction with the anticommutator spin algebra $\mathcal{A}$.  The algebra $\mathcal{A}$ is the unital associative $\K$-algebra defined by generators $x,y,z$ and relations\[xy+yx=2z,\qquad yz+zy=2x,\qquad zx+xz=2y.\]

Let $D$ denote a nonnegative integer, let $Q_{D}$ denote the hypercube of diameter $D$ and let $\tilde{Q}_{D}$ denote the antipodal quotient.  Let $T$ (resp. $\tilde{T}$) denote the Terwilliger algebra for $Q_{D}$ (resp. $\tilde{Q}_{D}$).

In this paper we obtain the following results.  When $D$ is even (resp. odd), we show that there exists a unique $\mathcal{A}$-module structure on the standard module for $Q_{D}$ (resp. $\tilde{Q}_{D}$) such that $x,y$ act as the adjacency and dual adjacency matrices respectively.  We classify the resulting irreducible $\mathcal{A}$-modules up to isomorphism.  We introduce weighted adjacency matrices  for $Q_{D}$ and $\tilde{Q}_{D}$.    When $D$ is even (resp. odd) we show that actions of the adjacency matrix, dual adjacency matrix and weighted adjacency matrix for $Q_{D}$ (resp. $\tilde{Q}_{D}$) on any irreducible $T$-module (resp. $\tilde{T}$-module) form a totally bipartite (resp. totally almost bipartite) Leonard triple of Bannai/Ito type and classify the Leonard triple up to isomorphism.
\end{abstract}

\section{Introduction}\label{S:intro}

Throughout this paper, $\K$ denotes an algebraically closed field of characteristic zero.

We now recall the definition of a Leonard triple.  To do this, we use the following terms.  A square matrix $B$ is said to be \emph{tridiagonal} whenever each nonzero entry lies on either the diagonal, the subdiagonal, or the superdiagonal.  Assume $B$ is tridiagonal.  Then $B$ is said to be \emph{irreducible} whenever each entry on the subdiagonal or superdiagonal is nonzero.

\begin{definition}\label{D:LT}
{\normalfont\cite[Definition 1.2]{curtin} Let $V$ denote a vector space over $\K$ with finite positive dimension.  By a \emph{Leonard triple} on $V$ we mean an ordered triple of linear transformations $A:V\to V$, $A^{*}:V\to V$, $A^{\varepsilon}:V\to V$ which satisfy the conditions (i)--(iii) below.
\begin{itemize}
\item[\rm (i)] There exists a basis for $V$ with respect to which the matrix representing $A$ is diagonal and the matrices representing $A^{*}$ and $A^{\varepsilon}$ are irreducible tridiagonal.
\item[\rm (ii)] There exists a basis for $V$ with respect to which the matrix representing $A^{*}$ is diagonal and the matrices representing $A^{\varepsilon}$ and $A$ are irreducible tridiagonal.
\item[\rm (iii)] There exists a basis for $V$ with respect to which the matrix representing $A^{\varepsilon}$ is diagonal and the matrices representing $A$ and $A^{*}$ are irreducible tridiagonal.
\end{itemize}
The \emph{diameter} of the Leonard triple $A,A^{*},A^{\varepsilon}$ is defined to be one less than the dimension of $V$.}
\end{definition}

Leonard triples were introduced by Curtin as a modification of Leonard pairs, introduced by Terwilliger (see \cite{terwilliger}).

We will be considering two families of Leonard pairs said to be totally bipartite and totally almost bipartite.  Before defining these families, we first review a few concepts.  Let $V$ denote a vector space over $\K$ with finite positive dimension.  By a \emph{decomposition} of $V$ we mean a sequence of one-dimensional subspaces of $V$ whose direct sum is $V$.  For any basis $\{v_{i}\}_{i=0}^{d}$ for $V$, the sequence $\{\K v_{i}\}_{i=0}^{d}$ is a decomposition of $V$; the decomposition $\{\K v_{i}\}_{i=0}^{d}$ is said to \emph{correspond} to the basis $\{v_{i}\}_{i=0}^{d}$.  Given a decomposition $\{V_{i}\}_{i=0}^{d}$ of $V$, for $0\leq i\leq d$ pick $0\ne v_{i}\in V_{i}$.  Then $\{v_{i}\}_{i=0}^{d}$ is a basis for $V$ which corresponds to $\{V_{i}\}_{i=0}^{d}$.

Let $A,A^{*},A^{\varepsilon}$ denote a Leonard triple on $V$.  A basis for $V$ is called \emph{standard} whenever it satisfies Definition \ref{D:LT}(i).  Observe that, given a decomposition $\{V_{i}\}_{i=0}^{d}$ of $V$, the following (i), (ii) are equivalent.
\begin{itemize}
\item[\rm (i)] There exists a standard basis for $V$ which corresponds to $\{V_{i}\}_{i=0}^{d}$.
\item[\rm (ii)] Every basis for $V$ which corresponds to $\{V_{i}\}_{i=0}^{d}$ is standard.
\end{itemize}
We say that the decomposition $\{V_{i}\}_{i=0}^{d}$ is \emph{standard} whenever (i), (ii) hold.  Observe that if the decomposition $\{V_{i}\}_{i=0}^{d}$ is standard, then so is $\{V_{d-i}\}_{i=0}^{d}$ and no other decomposition of $V$ is standard.

For any nonnegative integer $d$ let $\mathrm{Mat}_{d+1}(\K)$ denote the $\K$-algebra consisting of all $d+1$ by $d+1$ matrices that have entries in $\K$.  We index the rows and columns by $0,1,\ldots,d$.

Let $B\in\mathrm{Mat}_{d+1}(\K)$ be tridiagonal.  We say that $B$ is \emph{bipartite} whenever $B_{ii}=0$ for $0\leq i\leq d$.  Let $B\in\mathrm{Mat}_{d+1}(\K)$ be tridiagonal.  We say that $B$ is \emph{almost bipartite} whenever exactly one of $B_{0,0},B_{d,d}$ is nonzero and $B_{ii}=0$ for $1\leq i\leq d-1$.

\begin{definition}\label{D:bipartiteLT}
{\normalfont In Definition \ref{D:LT} we defined a Leonard triple $A,A^{*},A^{\varepsilon}$.  In that definition we mentioned six tridiagonal matrices.  The Leonard triple $A,A^{*},A^{\varepsilon}$ is said to be \emph{totally bipartite} (resp. \emph{totally almost bipartite}) whenever each of the six tridiagonal matrices is bipartite (resp. almost bipartite).}
\end{definition}

For notational convenience, we say that a Leonard triple is totally B/AB whenever it is either totally bipartite or totally almost bipartite.

In \cite{brown} we classified a family of totally B/AB Leonard triples said to have Bannai/Ito type.  We showed that totally B/AB Leonard triples of Bannai/Ito type satisfy the following property.

\begin{lemma}{\normalfont \cite[Theorem 9.5]{brown}}\label{L:nus}
Let $A,A^{*},A^{\varepsilon}$ denote a totally B/AB Leonard triple of Bannai/Ito type with diameter $d\geq3$.  Then there exist unique scalars $\nu,\nu^{*},\nu^{\varepsilon}\in\K$ such that
\begin{align}
AA^{*}+A^{*}A&=\nu^{\varepsilon}A^{\varepsilon}\label{E:nu1},\\
A^{*}A^{\varepsilon}+A^{\varepsilon}A^{*}&=\nu A\label{E:nu2},\\
A^{\varepsilon}A+AA^{\varepsilon}&=\nu^{*}A^{*}\label{E:nu3}.
\end{align}
\end{lemma}
From relations (\ref{E:nu1})--(\ref{E:nu3}) we showed a correspondence between the totally B/AB Leonard triples of Bannai/Ito type and the finite-dimensional irreducible representations of the anticommutator spin algebra $\mathcal{A}$.  The anticommutator spin algebra is a unital associative algebra over $\K$ with generators $x,y,z$ and relations
\begin{equation}\label{E:beginning}
xy+yx=2z,\qquad yz+zy=2x,\qquad zx+xz=2y.
\end{equation}
The algebra $\mathcal{A}$ was introduced by Arik and Kayserilioglu in \cite[Section 1]{arik-kayser}, in conjunction with fermionic quantum systems and the angular momentum algebra.

Other families of totally B/AB Leonard triples appear in the literature.  In \cite{HKP}, Havl\'{i}\v{c}ek, Klimyk and Po\v{s}ta displayed representations of the nonstandard $q$-deformed cyclically symmetric algebra $U'_{q}(\mathfrak{so}_{3})$.  These representations yield both totally bipartite and totally almost bipartite Leonard triples.  These Leonard triples are said to be of $q$-Racah type.  In \cite{miklavic}, Miklavi\v{c} studied totally bipartite Leonard triples associated with some representations of the Lie algebra $\mathfrak{sl}_{2}$ constructed using hypercubes.  These Leonard triples are said to be of Krawtchouk type.

The present paper is about how totally B/AB Leonard triples and finite-dimensional irreducible $\mathcal{A}$-modules arise naturally in the context of distance-regular graphs.   In \cite{go}, Go displayed an action of the Lie algebra $\mathfrak{sl}_{2}$ on the standard modules for hypercubes.  In \cite{miklavic}, Miklavi\v{c} used Go's results to display examples totally bipartite Leonard triples of Krawtchouk type.  We introduce an operator, called a skew operator, that gives $\mathcal{A}$-module structures to finite-dimensional $\mathfrak{sl}_{2}$-modules.  Using skew operators, we display an action of the algebra $\mathcal{A}$ on the primary modules of even-diameter hypercubes and the antipodal quotients of odd-diameter hypercubes similar to the $\mathfrak{sl}_{2}$-action from \cite{go}.  From the even-diameter hypercubes we obtain totally bipartite Leonard triples of Bannai/Ito type.  From the antipodal quotients of odd-diameter hypercubes we obtain totally almost bipartite Leonard triples of Bannai/Ito type.

Our main results will use the following notation.  Let $D$ denote a nonnegative integer.  Let $Q_{D}$ denote the hypercube of diameter $D$ and let $\tilde{Q}_{D}$ denote the antipodal quotient of $Q_{D}$.  Let $V$ (resp. $\tilde{V}$) denote the standard module for $Q_{D}$ (resp. $\tilde{Q}_{D}$) and Let $T$ (resp. $\tilde{T}$) denote the Terwilliger algebra of $Q_{D}$ (resp. $\tilde{Q}_{D}$)  When $D$ is even, let $\mathbf{A},\mathbf{B}$ denote the adjacency and dual adjacency matrices respectively for $Q_{D}$.  When $D$ is odd, let $\tilde{\mathbf{A}},\tilde{\mathbf{B}}$ denote the adjacency and dual adjacency matrices respectively for $\tilde{Q}_{D}$.

Our main results are as follows.  When $D$ is even, we show that there is a unique $\mathcal{A}$-module structure on $V$ such that the generators $x,y$ act as $\mathbf{A},\mathbf{B}$ respectively.  Let $W$ denote an irreducible $T$-module.  We show that $W$ is an irreducible $\mathcal{A}$-module classify the $\mathcal{A}$-module up to isomorphism.  When $D$ is odd, we show that there is a unique $\mathcal{A}$-module structure on $\tilde{V}$ such that the generators $x,y$ act as $\tilde{\mathbf{A}},\tilde{\mathbf{B}}$ respectively.  Let $\tilde{W}$ denote an irreducible $\tilde{T}$-module.  We show that $\tilde{W}$ is an irreducible $\mathcal{A}$-module classify the $\mathcal{A}$-module up to isomorphism.  We define weighted adjacency matrices $\mathbf{C},\tilde{\mathbf{C}}$ for $Q_{D},\tilde{Q}_{D}$ respectively.  When $D$ is even, we show that the $\mathbf{A},\mathbf{B},\mathbf{C}$ act on  $W$ as a totally bipartite Leonard triple of Bannai/Ito type and classify the Leonard triple up to isomorphism.  When $D$ is odd, we show that the $\tilde{\mathbf{A}},\tilde{\mathbf{B}},\tilde{\mathbf{C}}$ act on  $\tilde{W}$ as a totally almost bipartite Leonard triple of Bannai/Ito type and classify the Leonard triple up to isomorphism.

The paper is organized as follows.  In Section \ref{S:acsa} we recall the anticommutator spin algebra, its finite-dimensional irreducible representations and their relationship with totally B/AB Leonard triples of Bannai/Ito type.  In Section \ref{S:normalization} we introduce the notion of normalized Leonard triples and classify the normalized totally B/AB Leonard triples of Bannai/Ito type.  In Section \ref{S:sl2} we recall the Lie algebra $\mathfrak{sl}_{2}(\K)$ and some useful results about the algebra.  In Section \ref{S:s} we introduce the notion of a skew operator and use skew operators to give $\mathcal{A}$-module structures to finite-dimensional $\mathfrak{sl}_{2}(\K)$-modules.  In Section \ref{S:drg} we recall the notion of distance-regular graphs and some of their properties.  In Section \ref{S:modules} we recall the Terwilliger algebra of a distance-regular graph.  In Section \ref{S:hypercube} we recall some properties of hypercubes.  In Section \ref{S:sl2cube} we recall the results of Go concerning hypercubles and $\mathfrak{sl}_{2}(\K)$-modules.  In Section \ref{S:qhypercube} we recall some properties of the antipodal quotients of hypercubes.  In Section \ref{S:modcubes} we describe the $\mathcal{A}$-modules arising from hypercubes.  In Section \ref{S:qmodcubes} we describe the $\mathcal{A}$-modules arising from the antipodal quotients of odd-diameter hypercubes.  In Section \ref{S:LTs} we describe the Leonard triples arising from hypercubes and their antipodal quotients.

%%%%%%%%%%%%%%%%%%%%%%%%%%%%%%%%%%%%%%%%%%%%%%%%%%%%%%%%%%
\section{The anticommutator spin algebra}\label{S:acsa}

In this section we recall the anticommutator spin algebra $\mathcal{A}$ and the relationship between finite-dimensional irreducible $\mathcal{A}$-modules and totally B/AB Leonard triples of Bannai/Ito type.

We now recall the algebra $\mathcal{A}$.

\begin{definition}{\normalfont\cite[Section 1]{arik-kayser}}\label{D:A}
{\normalfont Let $\mathcal{A}$ denote the unital associative algebra over $\K$ with generators $x,y,z$ and relations
\begin{align}
xy+yx&=2z,\label{E:rel1}\\
yz+zy&=2x,\label{E:rel2}\\
zx+xz&=2y.\label{E:rel3}
\end{align}
We refer to the algebra $\mathcal{A}$ as the \emph{anticommutator spin algebra}.}
\end{definition}

The finite-dimensional irreducible $\mathcal{A}$-modules were classified up to isomorphism in \cite{brown}.  They fall into five families, described below.

\begin{lemma}{\normalfont \cite[Lemma 3.11]{brown}}\label{L:Bd}
Let $d$ denote a nonnegative even integer.  There exists an $\mathcal{A}$-module $V$ with basis $\{v_{i}\}_{i=0}^{d}$ on which $x,y,z$ act as follows.  For $0\leq i\leq d$,
\begin{align}
y.v_{i}=&(d-i+1)v_{i-1}+(i+1)v_{i+1},\label{E:Bd1}\\
x.v_{i}=&(-1)^{i}(d-2i)v_{i},\label{E:Bd2}\\
z.v_{i}=&(-1)^{i-1}(d-i+1)v_{i-1}+(-1)^{i}(i+1)v_{i+1},\label{E:Bd3}
\end{align}
where $v_{-1}=0$ and $v_{d+1}=0$.  The $\mathcal{A}$-module $V$ is irreducible.  An $\mathcal{A}$-module isomorphic to $V$ is said to have type $B(d)$.
\end{lemma}

\begin{lemma}{\normalfont \cite[Lemma 3.13]{brown}}\label{L:ABd0}
Let $d$ denote a nonnegative integer.  There exists an $\mathcal{A}$-module $V$ with basis $\{v_{i}\}_{i=0}^{d}$ on which $x,y,z$ act as follows.  For $0\leq i\leq d$,
\begin{align}
x.v_{i}=&(-1)^{d}(2d-i+2)v_{i-1}+(-1)^{d}(i+1)v_{i+1},\label{E:ABd01}\\
y.v_{i}=&(-1)^{d+i}(2d-2i+1)v_{i},\label{E:ABd02}\\
z.v_{i}=&(-1)^{i-1}(2d-i+2)v_{i-1}+(-1)^{i}(i+1)v_{i+1},\label{E:ABd03}
\end{align}
where $v_{-1}=0$ and $v_{d+1}=v_{d}$.  The $\mathcal{A}$-module $V$ is irreducible.  An $\mathcal{A}$-module isomorphic to $V$ is said to have type $AB(d,0)$.
\end{lemma}

\begin{lemma}{\normalfont \cite[Lemma 3.14]{brown}}\label{L:ABdx}
Let $d$ denote a nonnegative integer.  There exists an $\mathcal{A}$-module $V$ with basis $\{v_{i}\}_{i=0}^{d}$ on which $x,y,z$ act as follows.  For $0\leq i\leq d$,
\begin{align}
x.v_{i}=&(-1)^{d}(2d-i+2)v_{i-1}+(-1)^{d}(i+1)v_{i+1},\label{E:ABdx1}\\
y.v_{i}=&(-1)^{d+i+1}(2d-2i+1)v_{i},\label{E:ABdx2}\\
z.v_{i}=&(-1)^{i}(2d-i+2)v_{i-1}+(-1)^{i+1}(i+1)v_{i+1},\label{E:ABdx3}
\end{align}
where $v_{-1}=0$ and $v_{d+1}=v_{d}$.  The $\mathcal{A}$-module $V$ is irreducible.  An $\mathcal{A}$-module isomorphic to $V$ is said to have type $AB(d,x)$.
\end{lemma}

\begin{lemma}{\normalfont \cite[Lemma 3.15]{brown}}\label{L:ABdy}
Let $d$ denote a nonnegative integer.  There exists an $\mathcal{A}$-module $V$ with basis $\{v_{i}\}_{i=0}^{d}$ on which $x,y,z$ act as follows.  For $0\leq i\leq d$,
\begin{align}
x.v_{i}=&(-1)^{d+1}(2d-i+2)v_{i-1}+(-1)^{d+1}(i+1)v_{i+1},\label{E:ABdy1}\\
y.v_{i}=&(-1)^{d+i}(2d-2i+1)v_{i},\label{E:ABdy2}\\
z.v_{i}=&(-1)^{i}(2d-i+2)v_{i-1}+(-1)^{i+1}(i+1)v_{i+1},\label{E:ABdy3}
\end{align}
where $v_{-1}=0$ and $v_{d+1}=v_{d}$.  The $\mathcal{A}$-module $V$ is irreducible.  An $\mathcal{A}$-module isomorphic to $V$ is said to have type $AB(d,y)$.
\end{lemma}

\begin{lemma}{\normalfont \cite[Lemma 3.16]{brown}}\label{L:ABdz}
Let $d$ denote a nonnegative integer.  There exists an $\mathcal{A}$-module $V$ with basis $\{v_{i}\}_{i=0}^{d}$ on which $x,y,z$ act as follows.  For $0\leq i\leq d$,
\begin{align}
x.v_{i}=&(-1)^{d+1}(2d-i+2)v_{i-1}+(-1)^{d+1}(i+1)v_{i+1},\label{E:ABdz1}\\
y.v_{i}=&(-1)^{d+i+1}(2d-2i+1)v_{i},\label{E:ABdz2}\\
z.v_{i}=&(-1)^{i-1}(2d-i+2)v_{i-1}+(-1)^{i}(i+1)v_{i+1},\label{E:ABdz3}
\end{align}
where $v_{-1}=0$ and $v_{d+1}=v_{d}$.  The $\mathcal{A}$-module $V$ is irreducible.  An $\mathcal{A}$-module isomorphic to $V$ is said to have type $AB(d,z)$.
\end{lemma}

\begin{lemma}{\normalfont \cite[Theorem 3.20]{brown}}\label{L:class}
Every finite-dimensional irreducible $\mathcal{A}$-module is isomorphic to exactly one of the modules from Lemmas \ref{L:Bd}--\ref{L:ABdz}.
\end{lemma}

\begin{definition}{\normalfont \cite[Definition 3.17]{brown}}\label{D:diameter}
{\normalfont Let $V$ denote a finite-dimensional irreducible $\mathcal{A}$-module from Lemmas \ref{L:Bd}--\ref{L:ABdz}.  We define the \emph{diameter} of $V$ to be one less than the dimension of $V$.  Thus $\mathcal{A}$-modules of types $B(d)$ and $AB(d,n)$ have diameter $d$.}
\end{definition}

\begin{definition}{\normalfont \cite[Definition 3.18]{brown}}\label{D:BAB}
{\normalfont An $\mathcal{A}$-module $V$ is said to have \emph{type} $B$ when there exists an even integer $d\geq 0$ such that $V$ is of type $B(d)$.  The module is said to have \emph{type} $AB$ when there exists an integer $d\geq 0$ and $n\in\I$ such that $V$ is of type $AB(d,n)$.}
\end{definition}

Let $d\geq0$ denote an integer, let $V$ denote a vector space over $\K$ with dimension $d+1$ and $A,A^{*},A^{\varepsilon}$ denote a Leonard triple over $V$.  Let $\{v_{i}\}_{i=0}^{d}$ denote a standard basis for $A,A^{*},A^{\varepsilon}$ and let $\theta_{i}$ denote the eigenvalue for $A$ corresponding to $v_{i}$.  In \cite[Definition 9.2, Lemma 9.3]{brown} a Leonard triple is said to be of \emph{Bannai/Ito type} whenever the $(\theta_{i-2}-\theta_{i+1})/(\theta_{i-1}-\theta_{i})=-1$ for $2\leq i\leq d-1$.

We recall the classification of totally B/AB Leonard triples of Bannai/Ito type.  The following Lemma and Corollary classify the totally bipartite Leonard triples of Bannai/Ito type and the next Lemma and Corollary classify the totally almost bipartite Leonard triples of Bannai/Ito type.

\begin{lemma}\label{L:BLtriples}{\normalfont \cite[Theorem 9.6]{brown}}
Let $d$ denote an integer at least $3$ and let $\nu,\nu^{*},\nu^{\varepsilon}$ denote scalars in $\K$.  Then the following (i), (ii) are equivalent.
\begin{itemize}
\item[\rm (i)] There exists a totally bipartite Leonard triple $A,A^{*},A^{\varepsilon}$ of Bannai/Ito type with diameter $d$ that satisfies equations (\ref{E:nu1})--(\ref{E:nu3}).

\item[\rm (ii)] The integer $d$ is even and the scalars $\nu,\nu^{*},\nu^{\varepsilon}$ are nonzero.
\end{itemize}
Moreover, assume  (i), (ii) hold.  Then the Leonard triple $A,A^{*},A^{\varepsilon}$ is unique up to isomorphism.
\end{lemma}

As a consequence of Lemmas \ref{L:nus}, \ref{L:BLtriples}, one obtains the following Corollary.

\begin{corollary}\label{C:BLtriples}
There is a bijection between the set of isomorphism classes of totally bipartite Leonard triples of Bannai/Ito type with diameter at least $3$ and sequences $d,\nu,\nu^{*},\nu^{\varepsilon}$ of scalars such that $d\geq3$ is an even integer and $\nu,\nu^{*},\nu^{\varepsilon}$ are nonzero.
\end{corollary}

As we explain in \cite{brown}, the scalars $\nu,\nu^{*},\nu^{\varepsilon}$ are not suitable for the classification in the almost bipartite case.  We recall the classification of totally almost bipartite Leonard triples of Bannai/Ito type.

\begin{lemma}\label{L:ABLtriples}{\normalfont \cite[Theorem 9.7]{brown}}
Let $d$ denote an integer at least $3$ and let $\tau,\tau^{*},\tau^{\varepsilon}$ denote scalars in $\K$.  Then the following (i), (ii) are equivalent.
\begin{itemize}
\item[\rm (i)] There exists a totally almost bipartite Leonard triple $A,A^{*},A^{\varepsilon}$ of Bannai/Ito type with diameter $d$, $\mathrm{tr}(A)=\tau$, $\mathrm{tr}(A^{*})=\tau^{*}$ and $\mathrm{tr}(A^{\varepsilon})=\tau^{\varepsilon}$.

\item[\rm (ii)] The scalars $\tau,\tau^{*},\tau^{\varepsilon}$ are nonzero.
\end{itemize}
Moreover, assume  (i), (ii) hold.  Then the Leonard triple $A,A^{*},A^{\varepsilon}$ is unique up to isomorphism.
\end{lemma}

As a consequence of Lemma \ref{L:ABLtriples} and the fact that every linear transformation has a unique trace, one obtains the following Corollary.

\begin{corollary}\label{C:ABLtriples}
There is a bijection between the set of isomorphism classes of totally almost bipartite Leonard triples of Bannai/Ito type with diameter at least $3$ and sequences $d,\tau,\tau^{*},\tau^{\varepsilon}$ of scalars such that $d\geq3$ is an integer and $\tau,\tau^{*},\tau^{\varepsilon}$ are nonzero.
\end{corollary}

We recall the following two Lemmas from \cite{brown}, which provide a correspondence between $\mathcal{A}$-modules and Leonard triples.

\begin{lemma}\label{L:modLTs}{\normalfont \cite[Corollary 7.6, Theorem 9.4]{brown}}
Let $V$ denote a finite-dimensional irreducible $\mathcal{A}$-module and let $\xi,\xi^{*},\xi^{\varepsilon}$ denote nonzero scalars in $\K$.  Let $A,A^{*},A^{\varepsilon}$ denote the actions on $V$ of $\xi x,\xi^{*}y,\xi^{\varepsilon}z$ respectively.  Then $A,A^{*},A^{\varepsilon}$ is a Leonard triple on $V$ of Bannai/Ito type.  If $V$ is of type $B$ then $A,A^{*},A^{\varepsilon}$ is totally bipartite.  If $V$ is of type $AB$ then $A,A^{*},A^{\varepsilon}$ is totally almost bipartite.
\end{lemma}

\begin{lemma}\label{L:Ltrips}{\normalfont \cite[Theorem 9.5]{brown}}
Let $V$ denote a vector space with finite dimension at least $4$ and let $A,A^{*},A^{\varepsilon}$ denote a totally B/AB Leonard triple on $V$ of Bannai/Ito type.  Then there exists an irreducible $\mathcal{A}$-module structure on $V$ and nonzero scalars $\xi,\xi^{*},\xi^{\varepsilon}$ in $\K$ such that $x,y,z$ act as $A\xi^{-1},A^{*}\xi^{*-1},A^{\varepsilon}\xi^{\varepsilon-1}$ respectively.  If $A,A^{*},A^{\varepsilon}$ is totally bipartite then the $\mathcal{A}$-module $V$ is of type $B$ and if $A,A^{*},A^{\varepsilon}$ is totally almost bipartite then the $\mathcal{A}$-module $V$ is of type $AB$.  The $\mathcal{A}$-module structure of $V$ is uniquely determined by the scalars $\xi,\xi^{*},\xi^{\varepsilon}$ and the scalars $\xi,\xi^{*},\xi^{\varepsilon}$ are unique up to changing the sign of two of them.
\end{lemma}

%%%%%%%%%%%%%%%%%%%%%%%%%%%%%%%%%%%%%%%%%%%%%%%%%%%%%%%%%%
\section{Normalized Leonard triples}\label{S:normalization}

Observe that, given a vector space $V$ over $\mathbb{K}$ of finite positive dimension, Lemmas \ref{L:modLTs}, \ref{L:Ltrips} provided a correspondence between $\mathcal{A}$-modules on $V$ and totally B/AB Leonard triples on $V$ of Bannai/Ito type.  However, this is not a bijection as multiple Leonard triples may correspond to the same $\mathcal{A}$-module.  In this section, we introduce the notion of a normalized Leonard triple which will allow us to obtain a bijection.

\begin{definition}\label{D:normalization}
{\normalfont Let $V$ denote a vector space over $\mathbb{K}$ with positive finite dimension and let $A,A^{*},A^{\varepsilon}$ denote a  totally B/AB Leonard triple on $V$ of Bannai/Ito type.  We say the Leonard triple $A,A^{*},A^{\varepsilon}$ is \emph{normalized} whenever the scalars $\nu,\nu^{*},\nu^{\varepsilon}$ from equations (\ref{E:nu1})--(\ref{E:nu3}) are all equal to $2$.
}\end{definition}

We now describe the relationship between finite-dimensional irreducible $\mathcal{A}$-modules and normalized totally B/AB Leonard triples of Bannai/Ito type.

\begin{proposition}\label{P:Anormalization1}
Let $d$ denote a nonnegative integer and let $V$ denote a finite-dimensional irreducible $\mathcal{A}$-module of diameter $d$.  Then the actions of $x,y,z$ on $V$ form a normalized totally B/AB Leonard triple of Bannai/Ito type with diameter $d$. If the $\mathcal{A}$-module $V$ is of type $B$, then the Leonard triple $x,y,z$ is bipartite. If the $\mathcal{A}$-module $V$ is of type $AB$, then the Leonard triple $x,y,z$ is almost bipartite.
\end{proposition}

\noindent {\it Proof:} 
Immediate.
\hfill $\Box$ \\

\begin{proposition}\label{P:Anormalization2}
Let $V$ denote a vector space over $\mathbb{K}$ with finite dimension $d+1\geq4$ and let  $A,A^{*},A^{\varepsilon}$ denote a normalized totally B/AB Leonard triple on $V$ of Bannai/Ito type.  Then there exists a unique $\mathcal{A}$-module structure on $V$ such that $x,y,z$ act on $V$ as $A,A^{*},A^{\varepsilon}$.  If $A,A^{*},A^{\varepsilon}$ is totally bipartite then $V$ is of type $B$.  If $A,A^{*},A^{\varepsilon}$ is totally almost bipartite then $V$ is of type $AB$.  Moreover, the Leonard triple $A,A^{*},A^{\varepsilon}$ and the $\mathcal{A}$-module both have diameter $d$.
\end{proposition}

\noindent {\it Proof:} 
Immediate.
\hfill $\Box$ \\

\begin{remark}\label{R:normalization}
Let $V$ denote a vector space over $\mathbb{K}$ with finite dimension at least $4$.  Combining Propositions \ref{P:Anormalization1} and \ref{P:Anormalization2}, we obtain a bijection between the set of irreducible $\mathcal{A}$-modules on $V$ and the set of normalized totally B/AB Leonard triples on $V$.  Under this bijection the $\mathcal{A}$-modules on $V$ of type $B$ correspond to the normalized totally bipartite Leonard triples on $V$ and the $\mathcal{A}$-modules on $V$ of type $AB$ correspond to the normalized totally almost bipartite Leonard triples on $V$.  Moreover, the bijection preserves diameter.
\end{remark}

We now describe how normalized totally B/AB Leonard triples of Bannai/Ito type relate to general totally B/AB Leonard triples of Bannai/Ito type.

\begin{lemma}\label{L:nuxi}
Let $V$ denote a vector space over $\mathbb{K}$ with finite dimension at least $4$ and let $A,A^{*},A^{\varepsilon}$ denote a totally B/AB Leonard triple on $V$ of Bannai/Ito type.  Let $\xi,\xi^{*},\xi^{\varepsilon}$ denote scalars in $\mathbb{K}$ and let $\nu,\nu^{*},\nu^{\varepsilon}$ denote the scalars from equations (\ref{E:nu1})--(\ref{E:nu3}).  Then the Leonard triple $\xi A,\xi^{*}A^{*},\xi^{\varepsilon}A^{\varepsilon}$ is normalized if and only if the following hold.
\[
\xi^{2}=4(\nu^{*}\nu^{\varepsilon})^{-1},\qquad\xi^{*2}=4(\nu^{\varepsilon}\nu)^{-1},\qquad\xi^{\varepsilon2}=4(\nu\nu^{*})^{-1},\qquad\xi\xi^{*}\xi^{\varepsilon}=8(\nu\nu^{*}\nu^{\varepsilon})^{-1}.
\]
\end{lemma}

\noindent {\it Proof:} 
Routine.
\hfill $\Box$ \\

\begin{proposition}\label{P:normalization}
Let $V$ denote a vector space over $\mathbb{K}$ with finite dimension at least $4$ and let $A,A^{*},A^{\varepsilon}$ denote a totally B/AB Leonard triple on $V$ of Bannai/Ito type.  Then there exist exactly four sequences of nonzero scalars $\xi,\xi^{*},\xi^{\varepsilon}\in\mathbb{K}$ such that the Leonard triple $\xi A,\xi^{*}A^{*},\xi^{\varepsilon}A^{\varepsilon}$ is normalized.  Whenever the Leonard triple $A,A^{*},A^{\varepsilon}$ is totally bipartite, the four resulting normalized Leonard triples are isomorphic.  Whenever the Leonard triple $A,A^{*},A^{\varepsilon}$ is totally almost bipartite, the four resulting normalized Leonard triples are mutually nonisomorphic.
\end{proposition}

\noindent {\it Proof:} 
The existence of exactly four sequences of scalars $\xi,\xi^{*},\xi^{\varepsilon}$ is a routine consequence of Lemma \ref{L:nuxi}.

Assume the Leonard triple $A,A^{*},A^{\varepsilon}$ is totally bipartite.  Then the four resulting normalized Leonard triples are totally bipartite and have the same diameter.  By this, Theorem \ref{T:BLtriples} and Definition \ref{D:normalization} one routinely finds that the four resulting normalized Leonard triples are isomorphic as desired.

Now assume the Leonard triple $A,A^{*},A^{\varepsilon}$ is totally almost bipartite.  By comparing the traces of the resulting normalized Leonard triples one routinely finds that they are mutually nonisomorphic as desired.
\hfill $\Box$ \\

Observe that each totally almost bipartite Leonard triple of Bannai/Ito type can be scaled to four mutually nonisomorphic normalized Leonard triples.  We now elaborate on the differences between these four Leonard triples.  Recall the set $\mathbb{I}=\{0,x,y,z\}$.

\begin{definition}\label{D:nnormalization}
{\normalfont Let $d$ denote a nonnegative integer, let $V$ denote a vector space over $\mathbb{K}$ of dimension $d+1$ and let $A,A^{*},A^{\varepsilon}$ denote a  totally almost bipartite Leonard triple on $V$ of Bannai/Ito type.  Let $n\in\mathbb{I}$.  We say the Leonard triple $A,A^{*},A^{\varepsilon}$ is \emph{$n$-normalized} whenever the traces of $A,A^{*},A^{\varepsilon}$ agree with the following table.
\begin{center}
\begin{tabular}{|c||c|c|c|}
\hline
 & $\mathrm{tr}(A)$ & $\mathrm{tr}(A^{*})$ & $\mathrm{tr}(A^{\varepsilon})$ \\ \hline\hline
 $n=0$ & $(-1)^{d}(d+1)$ & $(-1)^{d}(d+1)$ & $(-1)^{d}(d+1)$ \\ \hline
 $n=x$ & $(-1)^{d}(d+1)$ & $(-1)^{d+1}(d+1)$ & $(-1)^{d+1}(d+1)$ \\ \hline
 $n=y$ & $(-1)^{d+1}(d+1)$ & $(-1)^{d}(d+1)$ & $(-1)^{d+1}(d+1)$ \\ \hline
 $n=z$ & $(-1)^{d+1}(d+1)$ & $(-1)^{d+1}(d+1)$ & $(-1)^{d}(d+1)$ \\ \hline
\end{tabular}
\end{center}
}\end{definition}

\begin{lemma}\label{L:2normalizations}
Let $A,A^{*},A^{\varepsilon}$ denote a totally almost bipartite Leonard triple of Bannai/Ito type.  Then the Leonard triple $A,A^{*},A^{\varepsilon}$ is normalized if and only if there exists $n\in\mathbb{I}$ such that the Leonard triple $A,A^{*},A^{\varepsilon}$ is $n$-normalized.
\end{lemma}

\noindent {\it Proof:} 
Immediate from Lemmas \ref{L:ABd0}--\ref{L:ABdz} and Definition \ref{D:nnormalization}.
\hfill $\Box$ \\

\begin{proposition}\label{P:nnormalization}
Let $V$ denote a vector space over $\mathbb{K}$ with finite dimension at least $4$ and let $A,A^{*},A^{\varepsilon}$ denote a totally almost bipartite Leonard triple on $V$ of Bannai/Ito type.  Then, for each $n\in\mathbb{I}$, there exists a unique sequence of nonzero scalars $\xi,\xi^{*},\xi^{\varepsilon}\in\mathbb{K}$ such that the Leonard triple $\xi A,\xi^{*}A^{*},\xi^{\varepsilon}A^{\varepsilon}$ is $n$-normalized.  These constitute the four sequences $\xi,\xi^{*},\xi^{\varepsilon}$ from Proposition \ref{P:normalization}.
\end{proposition}

\noindent {\it Proof:} 
Immediate
\hfill $\Box$ \\

\begin{proposition}\label{P:nbijection}
Under the bijection from Remark \ref{R:normalization}, the $\mathcal{A}$-modules on $V$ of type $AB(d,n)$ correspond to the $n$-normalized totally almost bipartite Leonard triples on $V$.
\end{proposition}

\noindent {\it Proof:} 
Immediate.
\hfill $\Box$ \\

%%%%%%%%%%%%%%%%%%%%%%%%%%%%%%%%%%%%%%%%%%%%%%%%%%%%%%%%%%
\section{The Lie algebra $\mathfrak{sl}_{2}(\K)$}\label{S:sl2}

In \cite{go}, Go displayed a relationship between the hypercube $Q_{D}$ and the finite-dimensional modules for the Lie algebra $\mathfrak{sl}_{2}$.  In \cite{miklavic}, Miklavi\v{c} used this relationship to display certain totally bipartite Leonard triples of Kroutchouk type.  We will use this information to obtain similar results involving the algebra $\mathcal{A}$ and totally B/AB Leonard triples of Bannai/Ito type.  First we recall the Lie algebra $\mathfrak{sl}_{2}$.  We will use a version of $\mathfrak{sl}_{2}$ favored by certain physicists (See \cite{zeidler}).

\begin{definition}\label{D:sl2}{\normalfont\cite[Proposition 1.21]{zeidler}}
{\normalfont We define $\mathfrak{sl}_{2}(\K)$ to be the Lie algebra over $\K$ with basis $X,Y,Z$ and Lie bracket}
\begin{align}
[X,Y]&=2\mathbf{i}Z\label{E:sl21}\\
[Y,Z]&=2\mathbf{i}X\label{E:sl22}\\
[Z,X]&=2\mathbf{i}Y\label{E:sl23}
\end{align}
\end{definition}
The finite-dimensional irreducible $\mathfrak{sl}_{2}$-modules are well-known (see \cite[Theorem, p. 33]{humphries}).  Up to isomorphism, there is one irreducible $\mathfrak{sl}_{2}$-module for each positive finite dimension.  Miklavi\v{c} showed that, on a $(d+1)$-dimensional $\mathfrak{sl}_{2}$-module $V$, $X,Y,Z$ act as a bipartite Leonard triple of diameter $d$

\begin{lemma}{\normalfont \cite[Lemma 9.1, Theorems 10.1--10.4]{miklavic}}\label{L:2bases}
Let $d$ denote a nonnegative integer and let $V$ denote a finite-dimensional irreducible $\mathfrak{sl}_{2}$-module of dimension $d+1$.  Then the following hold:
\begin{enumerate}
\item[\rm (i)] There exists a basis $\{v_{i}\}_{i=0}^{d}$ for $V$ on which $X,Y,Z$ act as follows.  For $0\leq i\leq d$:
\begin{eqnarray}
X.v_{i}&=&(d-i+1)v_{i-1}+(i+1)v_{i+1}\label{E:v1}\\
Y.v_{i}&=&(d-2i)v_{i}\label{E:v2}\\
Z.v_{i}&=&\mathbf{i}(d-i+1)v_{i-1}-\mathbf{i}(i+1)v_{i+1}\label{E:v3}
\end{eqnarray}
where $v_{-1}=0$ and $v_{d+1}=0$.

\item[\rm (ii)] There exists a basis $\{w_{i}\}_{i=0}^{d}$ for $V$ on which $X,Y,Z$ act as follows.  For $0\leq i\leq d$:
\begin{eqnarray}
X.w_{i}&=&\mathbf{i}(d-i+1)w_{i-1}-\mathbf{i}(i+1)w_{i+1}\label{E:w3}\\
Y.w_{i}&=&(d-i+1)w_{i-1}+(i+1)w_{i+1}\label{E:w2}\\
Z.w_{i}&=&(d-2i)w_{i}\label{E:w1}
\end{eqnarray}
where $v_{-1}=0$ and $v_{d+1}=0$.
\end{enumerate}
\end{lemma}

The $\mathfrak{sl}_{2}(\K)$-module $V$ from Lemma \ref{L:2bases} is said to have \emph{diameter $d$}.

\begin{lemma}\label{L:sum}
Let $V$ denote an irreducible $\mathfrak{sl}_{2}(\K)$-module with diameter $d$ and let $\{w_{i}\}_{i=0}^{d}$ denote the basis from Lemma \ref{L:2bases}(ii).  Then  $X.\sum_{i=0}^{d}w_{i}=d\sum_{i=0}^{d}w_{i}$.
\end{lemma}

\noindent {\it Proof:} 
Routine.
\hfill $\Box$ \\

%%%%%%%%%%%%%%%%%%%%%%%%%%%%%%%%%%%%%%%%%%%%%%%%%%%%%%%%%%
\section{From $\mathfrak{sl}_{2}(\K)$-modules to $\mathcal{A}$-modules}\label{S:s}

In this section we will demonstrate a relationship between $\mathcal{A}$-modules and $\mathfrak{sl}_{2}(\K)$-modules.  We introduce a type of operator called a skew operator. We use skew operators to give $\mathcal{A}$-module structures to each finite-dimensional $\mathfrak{sl}_{2}(\K)$-module.

\begin{definition}\label{D:skew}
{\normalfont By a \emph{skew operator} we mean a $\mathbb{K}$-linear operator $\sigma$ that acts on finite-dimensional $\mathfrak{sl}_{2}$-modules and satisfies the following relations:

\begin{eqnarray}
\sigma^{2}&=&I,\label{E:srel1}\\
\sigma X=-X\sigma,\qquad\qquad \sigma Y&=&Y\sigma,\qquad\qquad \sigma Z=-Z\sigma.\label{E:srel2}
\end{eqnarray}
Here, $X,Y,Z$ are from Definition \ref{D:sl2}.}
\end{definition}

\begin{lemma}\label{L:s}
Let $\sigma$ be a skew operator.  Then $X,\sigma Y,-\sigma\mathbf{i}Z$ satisfy relations (\ref{E:rel1})--(\ref{E:rel3}).
\end{lemma}

\noindent {\it Proof:} 
Routine.
\hfill $\Box$ \\

\begin{lemma}\label{L:diique}
Let $\sigma$ be a skew operator.  Then $-\sigma$ is a skew operator.  Furthermore every skew operator acts as either $\sigma$ or $-\sigma$ on any given finite-dimensional irreducible $\mathfrak{sl}_{2}$-module.
\end{lemma}

\noindent {\it Proof:} 
Routine calculation shows that $-\sigma$ is a skew operator.  Now let $\check{\sigma}$ be a skew operator and let $V$ be a finite-dimensional irreducible $\mathfrak{sl}_{2}$-module.  By (\ref{E:srel2}) $\sigma\check{\sigma}$ commutes with each of $X,Y,Z$.  Since the $\mathfrak{sl}_{2}$-module $V$ is irreducible, there exists $c\in\mathbb{K}$ such that $\sigma\hat{s}=cI$ on $V$.  Therefore $\check{\sigma}=c\sigma$ in view of (\ref{E:srel1}).  Also by (\ref{E:srel1}) $(c\sigma)^{2}=\check{\sigma}^{2}=I$, so $c=\pm1$.
\hfill $\Box$ \\

We have been discussing skew operators.  However we have not yet shown that a skew operator exists.  We now show that a skew operator exists by constructing one.  To construct our operator, we define two other operators, called $h$ and $k$.  Our skew operator will be equal to $hk$.  We will make use of the following lemma.

\begin{lemma}\label{L:nilp}
Let $V$ be a finite-dimensional $\mathfrak{sl}_{2}$-module.  Then $\frac{\mathbf{i}Y-X}{2},\frac{\mathbf{i}Y+X}{2}$ are nilpotent on $V$.
\end{lemma}

\noindent {\it Proof:} 
Finite-dimensional $\mathfrak{sl}_{2}$-modules are completely reducible.  Therefore we may assume without loss of generality that $V$ is irreducible.  Let $d$ be the diameter of $V$ and let $\{w_{i}\}_{i=0}^{d}$ be the basis for $V$ from Lemma \ref{L:2bases}(ii).  Then the matrix representing $\frac{\mathbf{i}Y-X}{2}$ With regards to $\{w_{i}\}_{i=0}^{d}$ is upper triangular with all diagonal entries zero and the matrix representing $\frac{\mathbf{i}Y+X}{2}$ is lower triangular with zeroes on the diagonal.  The result follows.
\hfill $\Box$ \\

\begin{definition}\label{D:p}
{\normalfont Let $h$ denote the $\mathbb{K}$-linear operator that acts on finite-dimensional $\mathfrak{sl}_{2}$-modules as follows:}
\[
h=\exp\left(\frac{\mathbf{i}Y-X}{2}\right)\exp\left(\frac{\mathbf{i}Y+X}{2}\right)\exp\left(\frac{\mathbf{i}Y-X}{2}\right).
\]
\end{definition}

We recall a fact about Lie algebras.  Let $\mathfrak{g}$ a Lie algebra and let $r,s \in\mathfrak{g}$ where $r$ is nilpotent on all finite-dimensional $\mathfrak{g}$-modules.  Then $\exp(r)s\exp(r)^{-1}=\exp(r)s\exp(-r)=$exp ad$(r)s$.

\begin{lemma}\label{L:pUp}
Let $h$ be as in Definition \ref{D:p}.  The following holds for $\xi\in\mathfrak{sl}_{2}$.
\begin{equation}\label{E:pUp}
h\xi h^{-1}=\mbox{\emph{exp ad}}\left(\frac{\mathbf{i}Y-X}{2}\right)\mbox{\emph{exp ad}}\left(\frac{\mathbf{i}Y+X}{2}\right)\mbox{\emph{exp ad}}\left(\frac{\mathbf{i}Y-X}{2}\right)\xi.
\end{equation}
\end{lemma}

\noindent {\it Proof:} 
Routine using Definition \ref{D:p} and the above comment. 
\hfill $\Box$ \\

\begin{lemma}\label{L:expad}
The following (i), (ii) hold.
\begin{itemize}
\item[\rm (i)] The matrix representing the action of \emph{exp ad}$(\frac{\mathbf{i}X-Y}{2})$ on $\mathfrak{sl}_{2}$ With regards to the basis $X,Y,Z$ is 
\[
\left(\begin{array}{ccc}
1/2 & \mathbf{i}/2 & -1\\
\mathbf{i}/2 & 3/2 & \mathbf{i}\\
1 & -\mathbf{i} & 1
\end{array}\right).
\]

\item[\rm (ii)] The matrix representing the action of \emph{exp ad}$(\frac{\mathbf{i}X+Y}{2})$ on $\mathfrak{sl}_{2}$ With regards to the basis $X,Y,Z$ is 
\[
\left(\begin{array}{ccc}
1/2 & -\mathbf{i}/2 & -1\\
-\mathbf{i}/2 & 3/2 & -\mathbf{i}\\
1 & \mathbf{i} & 1
\end{array}\right).
\]
\end{itemize}
\end{lemma}

\noindent {\it Proof:} 
Routine calculation using equations (\ref{E:sl21})--(\ref{E:sl23}).
\hfill $\Box$ \\

\begin{lemma}\label{L:pXY}
Let $X,Y,Z$ be as in Definition \ref{D:sl2}.  Then the element $h$ from Definition \ref{D:p} satisfies the following relations.

\begin{equation}
hX=-Xh\qquad\qquad hY=Yh\qquad\qquad hZ=-Zh\label{E:prel2}
\end{equation}
\end{lemma}

\noindent {\it Proof:} 
Evaluate each of $hXh^{-1},hYh^{-1},hZh^{-1}$ using Lemmas \ref{L:pUp}, \ref{L:expad}.
\hfill $\Box$ \\

A consequence of equation (\ref{E:prel2}) is that $h$ satisfies equation (\ref{E:srel2}). We now consider $h^{2}$.  To compute $h^{2}$ we look at how $h$ acts on finite-dimensional $\mathfrak{sl}_{2}$-modules.  We now examine this action in detail.

\begin{lemma}\label{L:scale}
Let $d$ denote a nonnegative integer and let $V$ denote a $\mathfrak{sl}_{2}$-module of diameter $d$. Pick an integer $i$ with $0\leq i\leq d$.  Let $v$ denote an eigenvector for $X$ with eigenvalue $d-2i$.  Then $v$ is an eigenvector for $h$ with eigenvalue $(-1)^{i}\mathbf{i}^{d}$.
\end{lemma}

\noindent {\it Proof:} 
Recall the basis $\{v_{i}\}_{i=0}^{d}$ for $V$ from Lemma \ref{L:2bases}(i) and recall that, for $0\leq i\leq d$, $v_{i}$ is an eigenvector for $X$ with eigenvalue $d-2i$.  By Lemma \ref{L:pXY}, $h$ commutes with $Y$.  Therefore $\{v_{i}\}_{i=0}^{d}$ is an eigenbasis for $h$.  By Lemma \ref{L:2bases}(i) the action of $X$ on $V$ is irreducible bipartite tridiagonal With regards to the basis $\{v_{i}\}_{i=0}^{d}$.  By  Lemma \ref{L:pXY}, $hX=-Xh$.  Therefore there exists a scalar $c\in\mathbb{K}$ such that, for $0\leq i\leq d$, $v_{i}$ is an eigenvector for $h$ with eigenvalue $(-1)^{i}c$.

We will now show that $c=\mathbf{i}^{d}$.  Recall the basis $\{w_{i}\}_{i=0}^{d}$ for $V$ from Lemma \ref{L:2bases}(ii).  By Lemma \ref{L:sum} and the results from the previous paragraph, $h\sum_{i=0}^{d}w_{i}=c\sum_{i=0}^{d}w_{i}$.  Through routine calculation one finds that, With regards to the basis $\{w_{i}\}_{i=0}^{d}$, the matrices representing $\exp(\frac{\mathbf{i}Y-X}{2}),\exp(\frac{\mathbf{i}Y+X}{2})$ are upper triangular.  Moreover, for $0\leq i\leq j\leq d$ the $(i,j)$-entry of the matrix representing $\exp(\frac{\mathbf{i}Y-X}{2})$ is ${d-i \choose d-j}\mathbf{i}^{j-i}$ and the $(i,j)$-entry of the matrix representing $\exp(\frac{\mathbf{i}Y+X}{2})$ is ${i \choose d-j}\mathbf{i}^{i-j}$.  From this we conclude, for $0\leq j\leq d$ the $(d,j)$-entry of the matrix representing $h$ with regards to the basis $\{w_{i}\}_{i=0}^{d}$ is
\begin{eqnarray*}
\sum_{k=0}^{j}\mathbf{i}^{d-k}{d \choose d-k}\mathbf{i}^{j-k}{d-k \choose d-j}&=&\mathbf{i}^{d+j}\sum_{k=0}^{j}(-1)^{2k}{d \choose d-k}{d-k \choose d-j}\\
&=&\mathbf{i}^{d+j}\sum_{k=0}^{j}(-1)^{2k}\frac{d!(d-k)!}{k!(d-k)!(d-j)!(j-k)!}\\
&=&\mathbf{i}^{d+j}\sum_{k=0}^{j}(-1)^{2k}\frac{d!j!}{k!j!(d-j)!(j-k)!}\\
&=&\mathbf{i}^{d+j}{d \choose j}\sum_{k=0}^{j}(-1)^{2k}{j \choose k}\\
&=&\delta_{j0}\mathbf{i}^{d}
\end{eqnarray*}
Therefore, when $h\sum_{i=0}^{d}v_{i}$ is expressed as a linear combination of $\{w_{i}\}_{i=0}^{d}$, the coefficient of $w_{d}$ must be $\mathbf{i}^{d}$, so $c=\mathbf{i}^{d}$.  The result follows.
\hfill $\Box$ \\

\begin{corollary}\label{C:h2central}
Let $h$ be from Definition \ref{D:p}.  Let $d$ denote a nonnegative integer and let $V$ denote a $\mathfrak{sl}_{2}$-module of diameter $d$.  Then $h^{2}=(-1)^{d}$ holds on $V$.
\end{corollary}

\noindent {\it Proof:}
Immediate.
\hfill $\Box$ \\

We see from Corollary \ref{C:h2central} that, although $h$ satisfies equation (\ref{E:srel2}), it does not satisfy equation (\ref{E:srel1}) and is therefore not a skew operator.  To construct our skew operator, we introduce the operator $k$.

\begin{lemma}\label{L:k}
There exists a unique $\mathbb{K}$-linear operator $k$ acting on finite-dimensional $\mathfrak{sl}_{2}$-modules such that $k$ acts as $1$ on each odd-dimensional irreducible submodule and as $-\mathbf{i}$ on each even-dimensional irreducible submodule.
\end{lemma}

\noindent {\it Proof:}
Follows from the fact that every finite-dimensional $\mathfrak{sl}_{2}$-module is completely irreducible.
\hfill $\Box$ \\

\begin{lemma}\label{L:kcentral}
Let $k$ denote the operator from Lemma \ref{L:k}.  Then $k$ commutes with each of $X,Y,Z$.
\end{lemma}

\noindent {\it Proof:}
Immediate.
\hfill $\Box$ \\

We are now ready to define the operator $s$.

\begin{definition}\label{D:s}
{\normalfont We define $s=hk$ where $h$ is from Definition \ref{D:p} and $k$ is from Lemma \ref{L:k}.}
\end{definition}

\begin{theorem}\label{T:sworks}
The operator $s$ from Definition \ref{D:s} is a skew operator.
\end{theorem}

\noindent {\it Proof:}
Equation (\ref{E:srel1}) is satisfied by Corollary \ref{C:h2central} and Lemma \ref{L:k}.  Equations (\ref{E:srel2}) are satisfied by equations (\ref{E:prel2}) and Lemma \ref{L:k}.
\hfill $\Box$ \\

\begin{theorem}\label{T:spk}
Let $V$ be a finite-dimensional $\mathfrak{sl}_{2}$-module.  Let $s$ be as in Definition \ref{D:s}.  Then the following (i), (ii) hold.

\begin{itemize}
\item[\rm (i)] There exists an $\mathcal{A}$-module structure on $V$ for which $x,y,z$ act as $X,sY-,s\mathbf{i}Z$ respectively.

\item[\rm (ii)] There exists an $\mathcal{A}$-module structure on $V$ for which $x,y,z$ act as $X,-sY,s\mathbf{i}Z$ respectively.
\end{itemize}
\end{theorem}

\noindent {\it Proof:}
(i) Immediate from Lemma \ref{L:s} and Theorem \ref{T:sworks}.

(ii) Immediate from Lemmas \ref{L:s}, \ref{L:diique} and Theorem \ref{T:sworks}.
\hfill $\Box$ \\

By the \emph{first} (resp. \emph{second}) $\mathcal{A}$-module structure on $V$ we mean the $\mathcal{A}$-module structure on $V$ from Theorem \ref{T:spk}(i) (resp. Theorem \ref{T:spk}(ii)).

Theorem \ref{T:spk} gives $\mathcal{A}$-module structures to every finite-dimensional $\mathfrak{sl}_{2}$-module.  When the $\mathfrak{sl}_{2}$-module is irreducible, it is not guaranteed that the corresponding $\mathcal{A}$-modules will be irreducible.  We now give necessary and sufficient conditions for the $\mathcal{A}$-module to be irreducible. When the $\mathcal{A}$-module isn't irreducible, we will describe how it decomposes into irreducible $\mathcal{A}$-modules. Moreover, we will show that every finite-dimensional irreducible $\mathcal{A}$-module can be obtained from an $\mathfrak{sl}_{2}$-module by this procedure.

\begin{theorem}\label{T:s-even}
Let $d$ denote a nonnegative even integer and let $V$ be a $(d+1)$-dimensional irreducible $\mathfrak{sl}_{2}$-module.  Then both $\mathcal{A}$-module structures on $V$ are irreducible and of type $B(d)$.
\end{theorem}

\noindent {\it Proof:}
Let $\{v_{i}\}_{i=0}^{d}$ be the basis for $V$ from Lemma \ref{L:2bases}(i).  For the first $\mathcal{A}$-module structure, observe the action of $X,sY,-s\mathbf{i}Z$ on $\{v_{i}\}_{i=0}^{d}$ coincides with the action of $x,y,z$ on the basis from Lemma \ref{L:Bd}. For the second $\mathcal{A}$-module structure, observe the action of $X,-sY,s\mathbf{i}Z$ on $\{v_{d-i}\}_{i=0}^{d}$ coincides with the action of $x,y,z$ on the basis from Lemma \ref{L:Bd}.
\hfill $\Box$ \\

\begin{theorem}\label{T:s-odd}
Let $d=2\delta+1$ denote a nonnegative odd integer and let $V$ be an irreducible $\mathfrak{sl}_{2}$-module with diameter $d$.  Then, under both $\mathcal{A}$-module structures, $V$ is the direct sum of the $\mathcal{A}$-submodules $V^{+}=\mathrm{span}\{v_{i}+v_{d-i}\}_{i=0}^{\delta}$ and $V^{-}=\mathrm{span}\{(-1)^{i}(v_{i}-v_{d-i})\}_{i=0}^{\delta}$ where $\{v_{i}\}_{i=0}^{d}$ is the basis from Lemma \ref{L:2bases}(i).  Under each $\mathcal{A}$-module structure, the $\mathcal{A}$-submodules $V^{+},V^{-}$ are irreducible.  Moreover, the following (i),(ii) hold.

\begin{itemize}
\item[\rm (i)] Under the first $\mathcal{A}$-module structure on $V$, the $\mathcal{A}$-submodule $V^{+}$ is of type $AB(\delta,0)$ and the $\mathcal{A}$-submodule $V^{-}$ is of type $AB(\delta,y)$.

\item[\rm (ii)] Under the second $\mathcal{A}$-module structure on $V$, the $\mathcal{A}$-submodule $V^{+}$ is of type $AB(\delta,z)$ and the $\mathcal{A}$-submodule $V^{-}$ is of type $AB(\delta,x)$.
\end{itemize}
\end{theorem}

\noindent {\it Proof:}
(i) The action of $sX,Y,s\mathbf{i}Z$ on $\{v_{i}+v_{d-i}\}_{i=0}^{\delta}$ coincides with the action of $x,y,z$ on the basis from Lemma \ref{L:ABd0} and the action of $X,sY,-s\mathbf{i}Z$ on $\{(-1)^{i}(v_{i}-v_{d-i})\}_{i=0}^{\delta}$ coincides with the action of $x,y,z$ on the basis from Lemma \ref{L:ABdy}.

(ii) The action of $-sX,Y,-s\mathbf{i}Z$ on $\{v_{i}+v_{d-i}\}_{i=0}^{\delta}$ coincides with the action of $x,y,z$ on the basis from Lemma \ref{L:ABdz} and the action of $X,-sY,s\mathbf{i}Z$ on $\{(-1)^{i}(v_{i}-v_{d-i})\}_{i=0}^{\delta}$ coincides with the action of $x,y,z$ on the basis from Lemma \ref{L:ABdx}.
\hfill $\Box$ \\

%%%%%%%%%%%%%%%%%%%%%%%%%%%%%%%%%%%%%%%%%%%%%%%%%%%%%%%%%%
\section{Distance-regular graphs}\label{S:drg}

In this section we review some definitions and basic results concerning distance-regular graphs.  For the remainder of the chapter we will work over the field $\mathbb{C}$.  Observe that, since $\mathbb{C}$ is an algebraically closed field of characteristic zero, every result from sections \ref{S:idempotents}--\ref{S:s} applies when $\mathbb{K}=\mathbb{C}$.

Let $\mathcal{X}$ denote a nonempty finite set.  Let $\mathrm{Mat}_{\mathcal{X}}(\mathbb{C})$ denote the $\mathbb{C}$-algebra of matrices with entries in $\mathbb{C}$ and with rows and columns indexed by $\mathcal{X}$.  Let $V=\mathbb{C}^{\mathcal{X}}$ denote the vector space over $\mathbb{C}$ consisting of column vectors with entries in $\mathbb{C}$ and rows indexed by $\mathcal{X}$.  We observe $\mathrm{Mat}_{\mathcal{X}}(\mathbb{C})$ acts on $V$ by left multiplication.  We refer to $V$ as the \emph{standard module} of $\mathrm{Mat}_{\mathcal{X}}(\mathbb{C})$.  For $\mathbf{x}\in\mathcal{X}$, we the vector in $V$ indexed by $\mathbf{x}$ is denoted $\hat{\mathbf{x}}$.  Let $L\in\mathrm{End}(V)$.  Given $\mathbf{N}\in\mathrm{Mat}_{\mathcal{X}}(\mathbb{C})$, we say $\mathbf{N}$ \emph{represents} $L$ whenever $L\hat{\mathbf{x}}=\sum_{\mathbf{y}\in\mathcal{X}}\mathbf{N}_{\mathbf{y}\mathbf{x}}\hat{\mathbf{y}}$ for all $\mathbf{x}\in\mathcal{X}$.

Let $\Gamma=(\mathcal{X},\mathcal{R})$ denote a finite, undirected, connected graph, without loops or multiple edges, with vertex set $\mathcal{X}$, edge set $\mathcal{R}$, path-length distance function $\partial$ and diameter $D=\mathrm{max}\{\partial(\mathbf{x},\mathbf{y})|\mathbf{x},\mathbf{y}\in \mathcal{X}\}$.  For a vertex $\mathbf{x}\in\mathcal{X}$ and an integer $i\geq0$ let $\Gamma_{i}(\mathbf{x})$ denote the set of vertices at distance $i$ from $a$.  For an integer $k\geq0$ we say $\Gamma$ is \emph{regular with valency $k$} whenever $|\Gamma_{1}(\mathbf{x})|=k$ for all $\mathbf{x}\in \mathcal{X}$.  We say $\Gamma$ is \emph{distance-regular} whenever for all integers $0\leq h,i,j\leq D$ and all $\mathbf{x},\mathbf{y}\in \mathcal{X}$ with $\partial(\mathbf{x},\mathbf{y})=h$ the number $p^{h}_{ij}$, defined to be $|\Gamma_{i}(\mathbf{x})\cap\Gamma_{j}(\mathbf{y})|$, is independent of $\mathbf{x},\mathbf{y}$.  The constants $p^{h}_{ij}$ are known as the \emph{intersection numbers} of $\Gamma$.  From now on we assume $\Gamma$ is distance-regular with $D\geq1$.  For convenience, set $c_{i}=p^{i}_{1,i-1}$ for $1\leq i\leq D$, $a_{i}=p^{i}_{1i}$ for $0\leq i\leq D$, $b_{i}=p^{i}_{1,i+1}$ for $0\leq i\leq D-1$, $k_{i}=p^{0}_{ii}$ for $0\leq i\leq D$, $c_{0}=0$ and $b_{D}=0$.  We observe that $\Gamma$ is regular with valency $k=k_{1}=b_{0}$ and that $c_{i}+a_{i}+b_{i}=k$ for $0\leq i\leq d$. By \cite[p. 127]{bcn} the following (i), (ii) hold for $0\leq h,i,j\leq D$.
\begin{enumerate}
\item[\rm (i)] $p^{h}_{ij}=0$ if one of $h,i,j$ is greater than the sum of the other two.

\item[\rm (ii)] $p^{h}_{ij}\ne0$ if one of $h,i,j$ is equal to the sum of the other two.
\end{enumerate}

We now recall the Bose-Mesner algebra of $\Gamma$.  For $0\leq i\leq D$ let $\mathbf{A}_{i}$ denote the matrix in $\mathrm{Mat}_{\mathcal{X}}(\mathbb{C})$ with entries
\[
(\mathbf{A}_{i})_{\mathbf{x}\mathbf{y}}=
\left\{\begin{array}{cl}
1 & \mbox{if $\partial(\mathbf{x},\mathbf{y})=i$,}\\
0 & \mbox{if $\partial(\mathbf{x},\mathbf{y})\ne i$,}
\end{array}\right.
\qquad(\mathbf{x},\mathbf{y}\in \mathcal{X}).
\]
We abbreviate $\mathbf{A}=\mathbf{A}_{1}$ and call this the \emph{adjacency matrix of $\Gamma$}.  Let $M$ denote the subalgebra of $\mathrm{Mat}_{\mathcal{X}}(\mathbb{C})$ generated by $\mathbf{A}$.  By \cite[p. 44]{bcn}, $\{\mathbf{A}_{i}\}_{i=0}^{D}$ is a basis for $M$.  We call $M$ the \emph{Bose-Mesner algebra of $\Gamma$}.  Observe that $M$ is commutative and semi-simple.  By \cite[p. 45]{bcn} there exists a basis $\{\mathbf{E}_{i}\}_{i=0}^{D}$ for $M$ such that
\begin{eqnarray}
\mathbf{E}_{0}&=&|\mathcal{X}|^{-1}\mathbf{J},\label{E:idem1}\\
\sum_{i=0}^{d}\mathbf{E}_{i}&=&\mathbf{I},\label{E:idem2}\\
\mathbf{E}_{i}\mathbf{E}_{j}&=&\delta_{ij}\qquad(0\leq i,j\leq D),\label{E:idem3}
\end{eqnarray}
where $\mathbf{I}$ and $\mathbf{J}$ denote the identity and the all-ones matrix of $\mathrm{Mat}_{\mathcal{X}}(\mathbb{C})$, respectively and where $\delta_{ij}$ denotes the Kronecker delta.  For convenience we define $E_{i}=0$ whenever $i<0$ or $i>d$.  The matrices $\{E_{i}\}_{i=0}^{D}$ are known as the \emph{primitive idempotents of $\Gamma$}, and $\mathbf{E}_{0}$ is called the \emph{trivial} idempotent.  We recall the eigenvalues of $\Gamma$.  Since $\{E_{i}\}_{i=0}^{D}$ is a basis for $M$, there exist scalars $\{\theta_{i}\}_{i=0}^{D}$ in $\mathbb{C}$ such that
\begin{equation}
\mathbf{A}=\sum_{i=0}^{D}\theta_{i}\mathbf{E}_{i}.
\end{equation}
Combining this with equations (\ref{E:idem1}) and (\ref{E:idem3}) we find that $\mathbf{A}\mathbf{E}_{i}=\mathbf{E}_{i}\mathbf{A}=\theta_{i}\mathbf{E}_{i}$ for $0\leq i\leq D$ and $\theta_{0}=k$.  Observe that $\{\theta_{i}\}_{i=0}^{D}$ are mutually distinct since $\mathbf{A}$ generates $M$.  We refer to $\theta_{i}$ as the \emph{eigenvalue of $\Gamma$ associated with $\mathbf{E}_{i}$}.  For $0\leq i\leq D$, let $m_{i}$ denote the rank of $\mathbf{E}_{i}$.  We call $m_{i}$ the \emph{multiplicity} of $\theta_{i}$.

By equations (\ref{E:idem2}), (\ref{E:idem3}),
\begin{equation}
V=\sum_{i=0}^{D}\mathbf{E}_{i}V\qquad(\mbox{direct sum}).
\end{equation}
By linear interpolation,
\begin{equation}\label{E:Es}
\mathbf{E}_{i}=\prod_{
\genfrac{}{}{0pt}{} {0\leq j\leq d}{j\ne i}
}\frac{\mathbf{A}-\theta_{j}\mathbf{I}}{\theta_{i}-\theta_{j}}\qquad(0\leq i\leq D).
\end{equation}
We now recall the $Q$-polynomial property.  Note that $\mathbf{A}_{i}\circ\mathbf{A}_{j}=\delta_{ij}\mathbf{A}_{i}$ for $0\leq i,j\leq D$, where $\circ$ is the entry-wise product.  Therefore $M$ is closed under $\circ$.  Thus there exist $q^{h}_{ij}\in\mathbb{C}$ for $0\leq h,i,j\leq D$ such that
\[
\mathbf{E}_{i}\circ\mathbf{E}_{j}=|\mathcal{X}|^{-1}\sum_{h=0}^{d}q^{h}_{ij}\mathbf{E}_{h}\qquad(0\leq i,j\leq D).
\]
The scalars $q^{h}_{ij}$ are called the \emph{Krein parameters} of $\Gamma$.  The ordering $\{\theta_{i}\}_{i=0}^{D}$ of eigenvalues for $\Gamma$ is said to be \emph{$Q$-polynomial} whenever the following (i), (ii) hold.
\begin{enumerate}
\item[\rm (i)] $q^{h}_{ij}=0$ if one of $h,i,j$ is greater than the sum of the other two.

\item[\rm (ii)] $q^{h}_{ij}\ne0$ if one of $h,i,j$ is equal to the sum of the other two.
\end{enumerate}
The graph $\Gamma$ is said to be \emph{$Q$-polynomial} whenever there exists a $Q$-polynomial ordering $\{\theta_{i}\}_{i=0}^{D}$ of the eigenvalues of $\Gamma$.

%%%%%%%%%%%%%%%%%%%%%%%%%%%%%%%%%%%%%%%%%%%%%%%%%%%%%%%%%%
\section{The Terwilliger algebra}\label{S:modules}

In this section we recall the dual Bose-Mesner algebra and the Terwilliger algebra of $\Gamma$.  For the rest of this section fix $\mathbf{x}\in\mathcal{X}$.  For $0\leq i\leq D$ let $\mathbf{E}_{i}^{*}=\mathbf{E}_{i}^{*}(\mathbf{x})$ denote the diagonal matrix in $\mathrm{Mat}_{\mathcal{X}}(\mathbb{C})$ with entries
\[
(\mathbf{E}_{i}^{*})_{\mathbf{y}\mathbf{y}}=
\left\{\begin{array}{cl}
1 & \mbox{if $\partial(\mathbf{x},\mathbf{y})=i$,}\\
0 & \mbox{if $\partial(\mathbf{x},\mathbf{y})\ne i$,}
\end{array}\right.
\qquad(\mathbf{y}\in \mathcal{X}).
\]
We call $\mathbf{E}_{i}^{*}$ the \emph{$i^{\mathrm{th}}$ dual idempotent of $\Gamma$ With regards to $\mathbf{x}$}.  We observe
\begin{eqnarray}
\sum_{i=0}^{d}\mathbf{E}^{*}_{i}&=&\mathbf{I},\label{E:didem2}\\
\mathbf{E}^{*}_{i}\mathbf{E}^{*}_{j}&=&\delta_{ij}\mathbf{E}^{*}_{i}\qquad(0\leq i,j\leq D),\label{E:didem3}
\end{eqnarray}
By construction $\{\mathbf{E}_{i}^{*}\}_{i=0}^{D}$ is linearly independent.  Let $M^{*}=M^{*}(\mathbf{x})$ denote the subalgebra of $\mathrm{Mat}_{\mathcal{X}}(\mathbb{C})$ spanned by $\{\mathbf{E}_{i}^{*}\}_{i=0}^{D}$.  We call $M^{*}$ the \emph{dual Bose-Mesner algebra of $\Gamma$ With regards to $\mathbf{x}$}.  We observe $M^{*}$ is commutative and semi-simple.

Assume $\Gamma$ is $Q$-polynomial with regards to the ordering $\{\theta_{i}\}_{i=0}^{D}$ of eigenvalues of $\Gamma$.  For $0\leq i\leq D$, let $\mathbf{A}_{i}^{*}=\mathbf{A}_{i}^{*}(\mathbf{x})$ denote the diagonal matrix in $\mathrm{Mat}_{\mathcal{X}}(\mathbb{C})$ with $(\mathbf{y},\mathbf{y})$-entry
\begin{equation}\label{E:Aistar}
(\mathbf{A}_{i}^{*})_{\mathbf{y}\mathbf{y}}=|\mathcal{X}|(\mathbf{E}_{i})_{\mathbf{x}\mathbf{y}}\qquad(\mathbf{y}\in\mathcal{X}).
\end{equation}
We call $\mathbf{A}^{*}$ the \emph{$i^{\mathrm{th}}$ dual distance matrix of $\Gamma$ With regards to $\mathbf{x}$}.  The matrix $\mathbf{A}_{1}^{*}$ is often denoted $\mathbf{A}^{*}$ and referred to as the \emph{dual adjacency matrix of $\Gamma$ with regards to $\mathbf{x}$}.  By \cite[Lemma 3.11]{bigterwilliger}, $M^{*}$ is generated by $\mathbf{A}^{*}$.  We recall the dual eigenvalues of $\Gamma$.  Since $\{\mathbf{E}_{i}^{*}\}_{i=0}^{D}$ is a basis for $M^{*}$ there exist scalars $\{\theta_{i}^{*}\}_{i=0}^{D}$ in $\mathbb{C}$ such that
\begin{equation}
\mathbf{A}^{*}=\sum_{i=0}^{D}\theta_{i}\mathbf{E}^{*}_{i}.
\end{equation}
Combining this with equation (\ref{E:didem3}) we find that $\mathbf{A}^{*}\mathbf{E}^{*}_{i}=\mathbf{E}^{*}_{i}\mathbf{A}^{*}=\theta^{*}_{i}\mathbf{E}^{*}_{i}$ for $0\leq i\leq D$.  The scalars are mutually distinct since $\mathbf{A}^{*}$ generates $M^{*}$.  Note that $\theta_{i}^{*}$ is an eigenvalue of $\mathbf{A}^{*}$ and $\mathbf{E}^{*}_{i}V$ is the corresponding eigenspace for $0\leq i\leq D$.  By equations (\ref{E:didem2}), (\ref{E:didem3}),
\begin{equation}
V=\sum_{i=0}^{D}\mathbf{E}^{*}_{i}V\qquad(\mbox{direct sum}).
\end{equation}
We call the sequence $\{\theta_{i}^{*}\}_{i=0}^{D}$ the \emph{dual eigenvalue sequence of $\Gamma$}.  Observe that for $0\leq i\leq D$ the rank of $\mathbf{E}^{*}_{i}$ is $k_{i}$.  Therefore $k_{i}$ is the multiplicity with which $\theta_{i}^{*}$ appears as an eigenvalue of $\mathbf{A}^{*}$.

By linear interpolation,
\begin{equation}\label{E:dEs}
\mathbf{E}^{*}_{i}=\prod_{
\genfrac{}{}{0pt}{} {0\leq j\leq D}{j\ne i}
}\frac{\mathbf{A}^{*}-\theta^{*}_{j}\mathbf{I}}{\theta^{*}_{i}-\theta^{*}_{j}}\qquad(0\leq i\leq d).
\end{equation}

By \cite[Lemma 3.2]{bigterwilliger} the following (i), (ii) hold for $0\leq h,j\leq D$.
\begin{enumerate}
\item[\rm (i)] $\mathbf{E}_{j}^{*}\mathbf{A}\mathbf{E}_{h}^{*}=0$ if and only if $p^{h}_{ij}=0$.

\item[\rm (ii)] $\mathbf{E}_{j}\mathbf{A}^{*}\mathbf{E}_{h}=0$ if and only if $q^{h}_{ij}=0$.
\end{enumerate}

Let $T=T(\mathbf{x})$ denote the subalgebra of $\mathrm{Mat}_{\mathcal{X}}(\mathbb{C})$ generated by $M$ and $M^{*}$.  We call $T$ the \emph{Terwilliger algebra of $\Gamma$ With regards to $x$}.  \cite[Definition 3.3]{bigterwilliger}

By a \emph{$T$-module} we mean a subspace $W$ of $V$ such that $\mathbf{B}W\subseteq W$ for all $\mathbf{B}\in T$.  Let $W$ denote a $T$-module.  Then $W$ is said to be \emph{irreducible} whenever $W$ is nonzero and $W$ contains no $T$-modules other than $0$ and $W$.

By \cite[Lemma 3.4(ii)]{bigterwilliger} $V$ decomposes into a direct sum of irreducible $T$-modules.  Let $W$ denote an irreducible $T$-module.  By  \cite[Lemma 3.4(iii)]{bigterwilliger} $W$ is the direct sum of the nonvanishing $\mathbf{E}_{i}W$ $(0\leq i\leq D)$ and the direct sum of the nonvanishing $\mathbf{E}^{*}_{i}V$ $(0\leq i\leq D)$.  By the \emph{endpoint} of $W$ we mean $\mathrm{min}\{i|0\leq i\leq D,\mathbf{E}_{i}^{*}W\ne0\}$.  By the \emph{diameter} of $W$ we mean $|\{i|0\leq i\leq D,\mathbf{E}_{i}^{*}W\ne0\}|-1$.  By the \emph{dual endpoint} of $W$ we mean $\mathrm{min}\{i|0\leq i\leq D,\mathbf{E}_{i}W\ne0\}$.  By the \emph{dual diameter} of $W$ we mean $|\{i|0\leq i\leq D,\mathbf{E}_{i}W\ne0\}|-1$.  By \cite[Lemma 4.5]{ITT} the diameter and the dual diameter of $W$ coincide.  Let $r$ and $r^{*}$ denote the endpoint and the dual endpoint of $W$, respectively, and let $d$ denote the diameter of $W$.  By  \cite[Lemma 3.9, Lemma 3.12]{bigterwilliger} the following (i), (ii) hold for $0\leq i\leq D$.
\begin{enumerate}
\item[\rm (i)] $\mathbf{E}_{i}W\ne0$ if and only if $r^{*}\leq i\leq r^{*}+d$.

\item[\rm (ii)] $\mathbf{E}^{*}_{i}W\ne0$ if and only if $r\leq i\leq r+d$.
\end{enumerate}

Let $W$ denote an irreducible $T$-module.  By  \cite[Lemma 3.9, Lemma 3.12]{bigterwilliger} the following (i), (ii) are equivalent.
\begin{enumerate}
\item[\rm (i)] $\mathrm{dim}(\mathbf{E}_{i}W)\leq1$ for $0\leq i\leq D$.

\item[\rm (ii)] $\mathrm{dim}(\mathbf{E}^{*}_{i}W)\leq1$ for $0\leq i\leq D$.
\end{enumerate}
In this case $W$ is called \emph{thin}.

%%%%%%%%%%%%%%%%%%%%%%%%%%%%%%%%%%%%%%%%%%%%%%%%%%%%%%%%%%
\section{The hypercube $Q_{D}$}\label{S:hypercube}

In this section we recall the hypercube graph and some of its basic properties.  Let $D$ denote a positive integer, and let $\{0,1\}^{D}$ denote the set of sequences $\{t_{i}\}_{i=1}^{D}$ there $t_{i}\in\{0,1\}$ for $1\leq i\leq D$.  Let $Q_{D}$ denote the graph with vertex set $\mathcal{X}=\{0,1\}^{D}$, and where two vertices are adjacent if and only if they differ in exactly one coordinate.  We call $Q_{D}$ the \emph{$D$-cube} or a \emph{hypercube}.  The graph $Q_{D}$ is connected and for $\mathbf{x},\mathbf{y}\in\mathcal{X}$ the distance $\partial(\mathbf{x},\mathbf{y})$ is the number of coordinates at which $\mathbf{x}$ and $\mathbf{y}$ differ.  In particular the diameter of $Q_{D}$ equals $D$.  The graph $Q_{D}$ is bipartite with bipartition $\mathcal{X}=\mathcal{X}^{+}\cup\mathcal{X}^{-}$, where $\mathcal{X}^{+}$ (resp. $\mathcal{X}^{-}$) is the set of vertices of $Q_{D}$ with an even (resp. odd) number of positive coordinates.

By \cite[p. 304]{bannai-ito} $Q_{D}$ is distance-regular with intersection numbers
\begin{equation}
a_{i}=0,\qquad b_{i}=D-i,\qquad c_{i}=i,\qquad k_{i}={D \choose i},\qquad(0\leq i\leq D).
\end{equation}
For $0\leq i\leq D$, $D-2i$ is an eigenvalue with multiplicity ${D \choose i}$.  The Graph $Q_{D}$ is $Q$-polynomial With regards to the ordering $\{D-2i\}_{i=0}^{D}$ of eigenvalues.  When $D$ is odd, this is the only $Q$-polynomial ordering of eigenvalues.  By \cite[p. 305]{bannai-ito} $D$ is even, $\{(-1)^{i}(D-2i)\}_{i=0}^{D}$ is also a $Q$-polynomial ordering of eigenvalues.

\begin{notation}\label{N:hypercube}
{\normalfont Let $D$ denote a nonnegative integer.  Let $Q_{D}$ denote the hypercube with diameter $D$.  Let $V$ denote the standard module for $Q_{D}$.  Let $M$ denote the Bose-Mesner algebra of $Q_{D}$.  Fix a vertex $\mathbf{x}$ of $Q_{D}$, let $M^{*}=M^{*}(\mathbf{x})$ denote the dual Bose-Mesner algebra of $Q_{D}$ With regards to $p$ and let $T=T(p)$ denote the Terwilliger algebra of $Q_{D}$ With regards to $\mathbf{x}$. For $0\leq i\leq D$, let $\mathbf{A}_{i}$ (resp. $\mathbf{A}_{i}^{*}$) denote the $i^{\mathrm{th}}$ distance matrix (resp. dual distance matrix) for $Q_{D}$ With regards to the $Q$-polynomial ordering $\{D-2i\}_{i=0}^{D}$ of eigenvalues with adjacency matrix $\mathbf{A}=\mathbf{A}_{1}$ and dual adjacency matrix $\mathbf{A}^{*}=\mathbf{A}^{*}_{1}$. For $0\leq i\leq D$, let $\mathbf{E}_{i}$ (resp. $\mathbf{E}_{i}^{*}$) denote the $i^{\mathrm{th}}$ primitive idempotent for $\mathbf{A}$ (resp. $\mathbf{A}^{*}$) With regards to the $Q$-polynomial ordering $\{D-2i\}_{i=0}^{D}$ of eigenvalues. For $0\leq i\leq D$, when $D$ is even, let $\mathbf{B}_{i}$ (resp. $\mathbf{F}_{i}$) denote the $i^{\mathrm{th}}$ dual distance matrix (resp. primitive idempotent of $\mathbf{A}$) associated with the $Q$-polynomial ordering $\{(-1)^{i}(D-2i)\}_{i=0}^{D}$ of eigenvalues with dual distance matrix $\mathbf{B}=\mathbf{B}_{1}$.}
\end{notation}

The following results about the hypercube $Q_{D}$ will be useful later.

\begin{lemma}\label{L:dminusi}
With regards to Notation \ref{N:hypercube}, the following (i), (ii) hold.
\begin{enumerate}
\item[\rm (i)] $(\mathbf{E}_{D-i})_{\mathbf{y}\mathbf{z}}=(-1)^{\partial(\mathbf{y},\mathbf{z})}(\mathbf{E}_{i})_{\mathbf{y}\mathbf{z}}\qquad(\mathbf{y},\mathbf{z}\in\mathcal{X})$.

\item[\rm (ii)] $\mathbf{A}_{D-i}^{*}$ is diagonal with $(\mathbf{A}_{D-i}^{*})_{\mathbf{y}\mathbf{y}}=(-1)^{\partial(\mathbf{x},\mathbf{y})}(\mathbf{A}_{i}^{*})_{\mathbf{y}\mathbf{y}}\qquad(\mathbf{y}\in\mathcal{X})$.
\end{enumerate}
\end{lemma}

\noindent {\it Proof:} 
(i) Let $\mathbf{S}\in\mathrm{Mat}_{\mathcal{X}}(\mathbb{C})$ denote the diagonal matrix whose $(\mathbf{y},\mathbf{y})$-entry is $1$ whenever $\mathbf{y}\in\mathcal{X}^{+}$ and $-1$ whenever $\mathbf{y}\in\mathcal{X}^{-}$.  Observe that the vector $v\in\mathbb{C}^{\mathcal{X}}$ is an eigenvector for $Q_{D}$ with eigenvalue $\theta$ if and only if $\mathbf{S}.v$ is an eigenvector for $Q_{D}$ with eigenvalue $-\theta$.  The result follows routinely from this and the fact that $\mathbf{E}_{i},\mathbf{E}_{D-i}$ project to the eigenspaces of $Q_{D}$ with eigenvalues $D-2i$ and $2i-D$ respectively.

(ii) Immediate from (i) and Equation (\ref{E:Aistar}).
\hfill $\Box$ \\

By \cite[(35)]{go}, the matrix $\mathbf{A}^{*}$ is diagonal with $(\mathbf{y},\mathbf{y})$-entry $D-2i$ where $i=\partial(\mathbf{x},\mathbf{y})$.  Combining this result with Lemma \ref{L:dminusi}(ii) we find that $\mathbf{A}_{D-1}^{*}$ is diagonal with $(\mathbf{y},\mathbf{y})$-entry $(-1)^{i}(D-2i)$ where $i=\partial(\mathbf{x},\mathbf{y})$.

%%%%%%%%%%%%%%%%%%%%%%%%%%%%%%%%%%%%%%%%%%%%%%%%%%%%%%%%%%
\section{$\mathfrak{sl}_{2}(\C)$-modules and hypercubes}\label{S:sl2cube}

In this section, we recall a result from Go relating hypercubes to $\mathfrak{sl}_{2}$-modules.

\begin{lemma}{\normalfont \cite[Theorem 13.2]{go}}\label{L:sl2cube}
With regard to Notation \ref{N:hypercube} there exists a unique $\mathfrak{sl}_{2}$-module structure on $V$ such that the generators $X,Y$ act as $\mathbf{A},\mathbf{A}^{*}$ respectively.
\end{lemma}

Observe that $X,Y$ generate $\mathfrak{sl}_{2}$ and $\mathbf{A},\mathbf{A}^{*}$ generate $T$.  Therefore the $\mathfrak{sl}_{2}$-module structure on $V$ from Lemma \ref{L:sl2cube} induces a surjective homomorphism of $\mathbb{C}$-algebras from $U(\mathfrak{sl}_{2})\to T$.  As a consequence we obtain the following Lemma.

\begin{lemma}\label{L:sl2cube2}
With regard to Notation \ref{N:hypercube} , let $V$ have the $\mathfrak{sl}_{2}$-module structure from Lemma \ref{L:sl2cube} and let $W\subseteq V$ denote an irreducible $T$-module.  Then $W$ is irreducible as an $\mathfrak{sl}_{2}$-module.
\end{lemma}

%%%%%%%%%%%%%%%%%%%%%%%%%%%%%%%%%%%%%%%%%%%%%%%%%%%%%%%%%%
\section{The antipodal quotient $\tilde{Q}_{D}$ of $Q_{D}$}\label{S:qhypercube}

In this section we recall the antipodal quotient of the hypercube and some of its basic properties.  The hypercube is an antipodal $2$-cover meaning that, for every vertex $\mathbf{y}\in\mathcal{X}$, there exists a unique vertex $\mathbf{y}'\in\mathcal{X}$ such that $\partial(\mathbf{y},\mathbf{y}')=D$.  The vertex $\mathbf{y}'$ is called the \emph{antipode} of $\mathbf{y}$.  Define the binary relation $\sim$ on $\mathcal{X}$ by $\mathbf{y}\sim\mathbf{z}$ whenever either $\mathbf{y}=\mathbf{z}$ or $\mathbf{y}'=\mathbf{z}$. Observe that $\sim$ is an equivalence relation. For every $\mathbf{y}\in\mathcal{X}$, the corresponding equivalence class will be denoted $[\mathbf{y}]$.  Let $\tilde{\mathcal{X}}$ denote the set of equivalence classes of $\sim$.  We define the graph $\tilde{Q}_{D}$ as follows.  The vertex set is $\mathcal{X}$.  Given $\mathbf{u},\mathbf{v}\in\tilde{\mathcal{X}}$, $\mathbf{u},\mathbf{v}$ are said to be adjacent in $\tilde{Q}_{D}$ if and only if there exist vertices $\mathbf{y},\mathbf{z}\in\mathcal{X}$ such that $\mathbf{y}\in\mathbf{u}$, $\mathbf{z}\in\mathbf{v}$ and $\mathbf{x},\mathbf{y}$ are adjacent in $Q_{D}$.

By  \cite[p. 306]{bannai-ito} $\tilde{Q}_{D}$ is distance-regular with diameter $\mathcal{D}$ where $\mathcal{D}=\lfloor\frac{D}{2}\rfloor$.  Observe that, when $D$ is even, $D=2\mathcal{D}$ and, when $D$ is odd, $D=2\mathcal{D}+1$. The intersection numbers of $\tilde{Q}_{D}$ are
\begin{eqnarray}
a_{i}=0,\qquad b_{i}=D-i,\qquad c_{i}=i,\qquad k_{i}={D \choose i},\qquad&(0\leq i\leq \mathcal{D}-1),\\
a_{\mathcal{D}}=0,\qquad b_{\mathcal{D}}=0,\qquad c_{\mathcal{D}}=D,\qquad k_{\mathcal{D}}=\frac{1}{2}{D \choose \mathcal{D}},\qquad&(\mbox{$D$ even}),\\
a_{\mathcal{D}}=\mathcal{\mathcal{D}}+1,\qquad b_{\mathcal{D}}=0,\qquad c_{\mathcal{D}}=\mathcal{D},\qquad k_{\mathcal{D}}={D \choose \mathcal{D}},\qquad&(\mbox{$D$ odd}).
\end{eqnarray}
For $0\leq i\leq\mathcal{D}$, $D-4i$ is an eigenvalue for $\tilde{Q}_{D}$ with multiplicity ${2D \choose2i}$.  The graph $\tilde{Q}_{D}$ is $Q$-polynomial With regards to the ordering $\{D-4i\}_{i=0}^{\mathcal{D}}$ of eigenvalues.  When $D$ is even, this is the only $Q$-polynomial ordering of eigenvalues.  By \cite[p. 305]{bannai-ito}, when $D$ is odd, $\{(-1)^{i}(D-2i)\}_{i=0}^{\mathcal{D}}$ is also a $Q$-polynomial ordering of the eigenvalues for $\tilde{Q}_{D}$.  This is the ordering that we will be concerned with in this thesis.

\begin{notation}\label{N:qhypercube}
{\normalfont Let the integer $D$, the graph $Q_{D}$ and the vertex $p$ be from Notation \ref{N:hypercube} and let $D=2\mathcal{D}+1$ be odd.  Let $\tilde{Q}_{D}$ denote the antipodal quotient of $Q_{D}$.  Let $\tilde{V}$ denote the standard module for $\tilde{Q}_{D}$.  Let $\tilde{M}$ denote the Bose-Mesner algebra of $\tilde{Q}_{D}$.  Let $\tilde{M}^{*}=\tilde{M}^{*}([\mathbf{x}])$ denote the dual Bose-Mesner algebra of $\tilde{Q}_{D}$ With regards to $[\mathbf{x}]$ and let $\tilde{T}=\tilde{T}([\mathbf{x}])$ denote the Terwilliger algebra of $\tilde{Q}_{D}$ With regards to $[\mathbf{x}]$. For $0\leq i\leq\mathcal{D}$, let $\tilde{\mathbf{A}}_{i}$ (resp. $\tilde{\mathbf{B}}_{i}$) denote the $i^{\mathrm{th}}$ distance matrix (resp. dual distance matrix) for $\tilde{Q}_{D}$ With regards to the $Q$-polynomial ordering $\{(-1)^{i}(D-2i)\}_{i=0}^{\mathcal{D}}$ of eigenvalues with adjacency matrix $\tilde{\mathbf{A}}=\tilde{\mathbf{A}}_{1}$ and dual adjacency matrix $\tilde{\mathbf{B}}=\tilde{\mathbf{B}}_{1}$. For $0\leq i\leq D$, let $\tilde{\mathbf{E}}_{i}$ (resp. $\tilde{\mathbf{F}}_{i}$) denote the $i^{\mathrm{th}}$ primitive idempotent for $\tilde{\mathbf{A}}$ (resp. $\tilde{\mathbf{B}}$) With regards to the $Q$-polynomial ordering $\{(-1)^{i}(D-2i)\}_{i=0}^{\mathcal{D}}$ of eigenvalues.}
\end{notation}

%%%%%%%%%%%%%%%%%%%%%%%%%%%%%%%%%%%%%%%%%%%%%%%%%%%%%%%%%%
\section{$\mathcal{A}$-modules and $Q_{D}$}\label{S:modcubes}

In Section \ref{S:s}, we displayed $\mathcal{A}$-module structures to finite-dimensional $\mathfrak{sl}_{2}$-modules. Go displayed $\mathfrak{sl}_{2}$-module structures to the primary module of the hypercube $Q_{D}$.  In this section we display an $\mathcal{A}$-module structure on the primary module of the hypercube $Q_{D}$.  We first recall the following result about irreducible $T$-modules.

\begin{lemma}{\normalfont \cite[Theorems 6.3, 8.1]{go}}\label{L:endpoint}
With regards to Notation \ref{N:hypercube}, let $W$ denote an irreducible $T$-module with endpoint $r$.  Then the following (i)--(iv) hold.
\begin{enumerate}
\item[\rm (i)] $r$ satisfies $0\leq r\leq D/2$.

\item[\rm (ii)] $W$ has diameter $D-2r$.

\item[\rm (iii)] $W$ had dual endpoint $r$.

\item[\rm (iii)] $W$ is thin.
\end{enumerate}
\end{lemma}

For the following Lemma, recall the operators $h$ from Definition \ref{D:p}, $k$ from Lemma \ref{L:k} and $s$ from Definition \ref{D:s}.

\begin{lemma}\label{L:psaction}
With regards to Notation \ref{N:hypercube}, let $V$ be endowed with the $\mathfrak{sl}_{2}$-module structure from Lemma \ref{L:sl2cube}.  Then the following (i)--(iii) hold.
\begin{itemize}
\item[\rm (i)] The matrix representing the action of $h$ on $V$ is diagonal with $(\mathbf{y},\mathbf{y})$-entry $(-1)^{i}\mathbf{i}^{D}$ for each vertex $\mathbf{y}$ of $Q_{D}$ where $i=\partial(\mathbf{x},\mathbf{y})$.

\item[\rm (ii)] The matrix representing the action of $k$ on $V$ is $I$ whenever $D$ is even and $-\mathbf{i}I$ whenever $D$ is odd.

\item[\rm (iii)] The matrix representing the action of $s$ on $V$ is diagonal with $(\mathbf{y},\mathbf{y})$-entry $(-1)^{\lfloor\frac{D}{2}\rfloor+i}$ for each vertex $\mathbf{y}$ of $Q_{D}$ where $i=\partial(\mathbf{x},\mathbf{y})$.
\end{itemize}
\end{lemma}

\noindent {\it Proof:} 
(i) By Lemma \ref{L:sl2cube}(i), the action of $Y$ on $V$ coincides with the action of $\mathbf{A}^{*}$ on $V$.  By the comment at the end of Section \ref{S:hypercube}, the matrix $\mathbf{A}^{*}$ is diagonal with $(\mathbf{y},\mathbf{y})$-entry $d-2i$.  The result follows from this information along with Lemma \ref{L:scale}.

(ii) By  Lemma \ref{L:endpoint}(ii), the dimension of every irreducible $T$-module has the same parity.  By this and Lemma \ref{L:k}, $k$ acts on $V$ as $1$ whenever $D$ is even and as $-\mathbf{i}$ whenever $D$ is odd, as desired.

(iii) Immediate from Definition \ref{D:s} and parts (i), (ii).
\hfill $\Box$ \\

\begin{lemma}\label{L:snsactions}
With regards to Notation \ref{N:hypercube}, let $V$ be endowed with the $\mathfrak{sl}_{2}$-module structure from Lemma \ref{L:sl2cube}.  Let $\varepsilon=1$ when $D$ is congruent to $0$ or $1$ modulo $4$ and $\varepsilon=-1$ when $D$ is congruent to $2$ or $3$ modulo $4$.  Then the following (i),(ii) hold.
\begin{itemize}
\item[\rm (i)] In the first $\mathcal{A}$-module structure on $V$ the generators $x,y$ act on $V$ as $\mathbf{A},\varepsilon\mathbf{A}_{D-1}^{*}$ respectively.

\item[\rm (ii)] In the second $\mathcal{A}$-module structure on $V$ the generators $x,y$ act on $V$ as $\mathbf{A},-\varepsilon\mathbf{A}_{D-1}^{*}$ respectively.
\end{itemize}
\end{lemma}

\noindent {\it Proof:} 
(i) By Lemma \ref{L:s}, the generators $x,y$ act on $V$ as $sX,Y$ respectively.  By Lemma \ref{L:sl2cube}(i), the generators $X,Y$ act as $\mathbf{A},\mathbf{A}^{*}$ respectively.  The matrix representing $s\mathbf{A}^{*}$ is diagonal with $(\mathbf{y},\mathbf{y})$-entry $(-1)^{\mathcal{D}+i}(D-2i)$ where $i=\partial(\mathbf{y},\mathbf{y})$ and $\mathcal{D}=\lfloor\frac{D}{2}\rfloor$.  By  the comment at the end of Section \ref{S:hypercube}, the matrix representing $\mathbf{A}^{*}_{D-1}$ is diagonal with $(\mathbf{y},\mathbf{y})$-entry $(-1)^{i}(D-2i)$ where $i=\partial(\mathbf{x},\mathbf{y})$.  The result follows.

(ii) Immediate from part (i) and Lemma \ref{L:s}.
\hfill $\Box$ \\

\begin{definition}\label{D:posnegmod}
{\normalfont With regards to Notation \ref{N:hypercube}, let $V$ be given an $\mathcal{A}$-module structure. We say the $\mathcal{A}$-module structure on $V$ is \emph{positive} whenever $x,y$ act on $V$ as $\mathbf{A},\mathbf{A}_{D-1}^{*}$ respectively. We say the $\mathcal{A}$-module structure on $V$ is \emph{negative} whenever $x,y$ act on $V$ as $\mathbf{A},-\mathbf{A}_{D-1}^{*}$ respectively}
\end{definition}

We now show the existence and uniqueness of positive and negative $\mathcal{A}$-module structures.

\begin{lemma}\label{L:ad-1}
With regards to Notation \ref{N:hypercube}, there exists a unique positive $\mathcal{A}$-module structure on $V$ and a unique negative $\mathcal{A}$-module structure on $V$.
\end{lemma}

\noindent {\it Proof:} 
By Lemma \ref{L:snsactions}, one finds examples of positive and negative $\mathcal{A}$-module structures on $V$ as given in the following table.
\begin{center}
\begin{tabular}{|c||c|c|}
\hline
 & $D$ congruent to $0$ or $1$, mod $4$ & $D$ congruent to $2$ or $3$ mod $4$ \\ \hline\hline
positive & first & second \\ \hline
negative & second & first \\ \hline
\end{tabular}
\end{center}
The $\mathcal{A}$-module structures are unique because $x,y$ generate $\mathcal{A}$.
\hfill $\Box$ \\

We will describe the irreducible $\mathcal{A}$-submodules of the positive $\mathcal{A}$-module structure on $V$.

\begin{theorem}\label{T:BILT}
With regards to notation \ref{N:hypercube}, assume $D$ is even.  The following (i), (ii) hold
\begin{itemize}
\item[\rm (i)] There is a unique $\mathcal{A}$-module structure on $V$ such that the generators $x,y$ act as $\mathbf{A},\mathbf{B}$ respectively.

\item[\rm (ii)] Let $W$ denote an irreducible $T$-module with endpoint $r$.  Then $W$ is an irreducible $\mathcal{A}$-module of type $B(D-2r)$.
\end{itemize}
\end{theorem}

\noindent {\it Proof:} 
(i) Immediate from Lemma \ref{L:ad-1}.

(ii) By  Lemma \ref{L:endpoint}(ii),(iv), $W$ has dimension $D-2r+1$.  The result follows from this information and Theorem \ref{T:s-even}.
\hfill $\Box$ \\

In Theorem \ref{T:oddcase1} we will describe the irreducible $\mathcal{A}$-modules when $D$ is odd.  To do this, we introduce some useful results.

\begin{lemma}{\normalfont \cite[Lemma 4.1]{collins}}\label{L:Ad1}
With regards to Notation \ref{N:hypercube}, the following (i), (ii) hold.
\begin{enumerate}
\item[\rm (i)] $\mathbf{A}_{D}\hat{\mathbf{y}}=\hat{\mathbf{y}}'\qquad(\mathbf{y}\in\mathcal{X})$.

\item[\rm (ii)] $\mathbf{A}_{D}:V\to V$ is an isomorphism of vector spaces.
\end{enumerate}
\end{lemma}

\begin{lemma}{\normalfont \cite[Lemma 4.2]{collins}}\label{L:Ad2}
With regards to Notation \ref{N:hypercube}, $\mathbf{A}_{D}^{2}=\mathbf{I}$.
\end{lemma}

\begin{lemma}{\normalfont \cite[Lemma 4.3]{collins}}\label{L:Vpm}
With regards to Notation \ref{N:hypercube}, define
\[
V_{+}=\mathrm{span}\{\hat{\mathbf{y}}+\hat{\mathbf{y}}'|\mathbf{y}\in\mathcal{X}\},\qquad\qquad V_{-}=\mathrm{span}\{\hat{\mathbf{y}}-\hat{\mathbf{y}}'|\mathbf{y}\in\mathcal{X}\}.
\]
Then the following (i)--(v) hold.
\begin{enumerate}
\item[\rm (i)] $V=V_{+}+V_{-}$ (direct sum).

\item[\rm (ii)] $(\mathbf{A}_{D}-\mathbf{I})V_{+}=0$.

\item[\rm (iii)] $(\mathbf{A}_{D}+\mathbf{I})V_{-}=0$.

\item[\rm (iv)] $(\mathbf{A}_{D}+\mathbf{I})V=V_{+}$.

\item[\rm (v)] $(\mathbf{A}_{D}-\mathbf{I})V=V_{-}$.
\end{enumerate}
\end{lemma}

\begin{lemma}{\normalfont \cite[Lemma 8.3]{collins}}\label{L:VpmE}
With regards to Notation \ref{N:hypercube}, let $V_{+},V_{-}$ be as in Lemma \ref{L:Vpm}.  Then the following (i), (ii) hold.
\begin{enumerate}
\item[\rm (i)] $V_{+}=\sum_{\mbox{$i$ \rm{even}}}\mathbf{E}_{i}V$ (direct sum).

\item[\rm (ii)] $V_{-}=\sum_{\mbox{$i$ \rm{odd}}}\mathbf{E}_{i}V$ (direct sum).
\end{enumerate}
\end{lemma}

\begin{theorem}\label{T:oddcase1}
With regards to Notation \ref{N:hypercube}, assume $D=2\mathcal{D}+1$ is odd.  Let $V$ be endowed with the $\mathcal{A}$-module structure from Lemma \ref{L:ad-1}.  Then the following (i), (ii) hold
\begin{itemize}
\item[\rm (i)] The subspaces $V_{+},V_{-}$ of $V$ from Lemma \ref{L:Vpm} are $\mathcal{A}$-submodules of $V$.

\item[\rm (ii)] Let $W$ denote an irreducible $T$-module with endpoint $r$.  Then $W$ is a direct sum of two irreducible $\mathcal{A}$-modules $(W\cap V_{+}),(W\cap V_{-})$.  The space $W\cap V_{+}$ is an irreducible $\mathcal{A}$-module of type $AB(\mathcal{D}-r,n)$ where $n$ is given in the following table.
\begin{center}
\begin{tabular}{|c||c|c|}
\hline
n & $\mathcal{D}$ even & $\mathcal{D}$ odd \\ \hline\hline
$r$ even & $0$ & $z$ \\ \hline
$r$ odd & $x$ & $y$ \\ \hline
\end{tabular}
\end{center}

Moreover $W\cap V_{-}$ is an irreducible $\mathcal{A}$-module of type $AB(\mathcal{D}-r,n)$ where $n$ is given in the following table.
\begin{center}
\begin{tabular}{|c||c|c|}
\hline
n & $\mathcal{D}$ even & $\mathcal{D}$ odd \\ \hline\hline
$r$ even & $y$ & $x$ \\ \hline
$r$ odd & $z$ & $0$ \\ \hline
\end{tabular}
\end{center}
\end{itemize}
\end{theorem}

\noindent {\it Proof:} 
(i)  One routinely finds that
\begin{eqnarray}
\mathbf{A}V_{+}\subseteq V_{+},&\qquad\qquad\mathbf{A}^{*}_{D-1}V_{+}\subseteq V_{+},\\
\mathbf{A}V_{+}\subseteq V_{-},&\qquad\qquad\mathbf{A}^{*}_{D-1}V_{+}\subseteq V_{-}.
\end{eqnarray}
The result follows from the above equations, Lemma \ref{L:ad-1} and the fact that $x,y$ generate $\mathcal{A}$.

(ii) By Lemma \ref{L:sl2cube2}, $W$ is irreducible as an $\mathfrak{sl}_{2}$-module and, by Lemma \ref{L:endpoint}(ii), $W$ has diameter $D-2r$.  By Theorem \ref{T:s-odd} and the fact that $D-2r$ is odd, one finds that $W$ is the direct sum of two irreducible $\mathcal{A}$-modules each with diameter $\mathcal{D}-r$.  By Theorem \ref{T:s-odd} and Lemma \ref{L:snsactions}, the two irreducible modules are of type $AB(\mathcal{D}-r,0)$ and $AB(\mathcal{D}-r,x)$ when $\mathcal{D}-r$ is even and they are of type $AB(\mathcal{D}-r,y)$ and $AB(\mathcal{D}-r,z)$ when $\mathcal{D}-r$ is odd.  By Lemmas \ref{L:VpmE}, \ref{L:endpoint}(iii), the $\mathcal{A}$-modules $W\cap V_{+}$ and $W\cap V_{-}$ each have diameter $\mathcal{D}-r$ and, by Lemma \ref{L:Vpm}, their direct sum is $W$.  Therefore $W\cap V_{+}$ and $W\cap V_{-}$ are the two irreducible $\mathcal{A}$-submodules of $W$ from Theorem \ref{T:s-odd}.  By Lemmas \ref{L:VpmE}, \ref{L:endpoint}(iii), the action of $A$ on $V_{+}$ has trace $(-1)^{r}(\mathcal{D}-r+1)$ and the action of $A$ on $V_{-}$ had trace $(-1)^{r+1}(\mathcal{D}-r+1)$.  Comparing these traces to the traces for $x$ in Lemmas \ref{L:ABd0}--\ref{L:ABdz}, one routinely obtains the result as desired.
\hfill $\Box$ \\

%%%%%%%%%%%%%%%%%%%%%%%%%%%%%%%%%%%%%%%%%%%%%%%%%%%%%%%%%%
\section{$\mathcal{A}$-modules and $\tilde{Q}_{D}$}\label{S:qmodcubes}

In Section \ref{S:modcubes}, we displayed an $\mathcal{A}$-module structure on the primary module of the hypercube $Q_{D}$.  In this section we will display an $\mathcal{A}$-module structure for the standard module of the antipodal quotient $\tilde{Q}_{D}$ and we will describe the irreducible $\mathcal{A}$-modules when $D$ is odd.  Before this we describe the irreducible $\tilde{T}$-modules.

\begin{definition}\label{D:psi}
{\normalfont With regards to Notation \ref{N:hypercube}, \ref{N:qhypercube}, let $\psi:V\to\tilde{V}$ denote the vector space homomorphism sending $\hat{\mathbf{y}}\mapsto\widehat{[\mathbf{y}]}$.}
\end{definition}

\begin{lemma}\label{L:psi}{\normalfont \cite[Lemma 10.4]{collins}}
With regards to Notation \ref{N:qhypercube}, let $\psi$ be from Definition \ref{D:psi}.  Then the following (i), (ii) hold.
\begin{itemize}
\item[\rm (i)]$V_{-}=\mathrm{Ker}\psi$.

\item[\rm (ii)] Let $\varphi$ denote the restriction of $\psi$ to $V_{+}$.  Then $\varphi:V_{+}\to\tilde{V}$ is an isomorphism of vector spaces.
\end{itemize}
\end{lemma}

\begin{lemma}{\normalfont \cite[Corollary 11.2]{collins}}\label{L:qbij}
With regards to Notation \ref{N:hypercube}, \ref{N:qhypercube}, the map $W\mapsto\psi(W)$ is a bijection from the set of $T$-modules to the set of $\tilde{T}$-modules.
\end{lemma}

\begin{lemma}\label{L:qendpoint}
With regards to Notation \ref{N:hypercube}, \ref{N:qhypercube}, let $W$ denote an irreducible $\tilde{T}$-module with endpoint $r$.  Then the following (i)--(iii) hold.
\begin{enumerate}
\item[\rm (i)] $W$ has diameter $\mathcal{D}-r$.

\item[\rm (ii)] $W$ has dual endpoint $r$.

\item[\rm (iii)] $W$ is thin.
\end{enumerate}
\end{lemma}

\noindent {\it Proof:} 
Follows from Lemmas \ref{L:endpoint}, \ref{L:VpmE}, \ref{L:psi}, \ref{L:qbij}.
\hfill $\Box$ \\

\begin{theorem}\label{T:oddcase2}
With regards to Notation \ref{N:qhypercube} assume $D=2\mathcal{D}+1$ is odd.  Then the following (i), (ii) hold.
\begin{itemize}
\item[\rm (i)] There is a unique $\mathcal{A}$-module structure on $\tilde{V}$ such that the generators $x,y$ act as $\tilde{\mathbf{A}},\tilde{\mathbf{B}}$ respectively.

\item[\rm (ii)] Let $W$ denote an irreducible $\tilde{T}$-module with endpoint $r$.  Then $W$ is an irreducible $\mathcal{A}$-module of type $AB(\mathcal{D}-r,n)$ where $n$ given in the following table.
\begin{center}
\begin{tabular}{|c||c|c|}
\hline
n & $\mathcal{D}$ even & $\mathcal{D}$ odd \\ \hline\hline
$r$ even & $0$ & $z$ \\ \hline
$r$ odd & $x$ & $y$ \\ \hline
\end{tabular}
\end{center}
\end{itemize}
\end{theorem}

\noindent {\it Proof:} 
(i) With regards to Notation \ref{N:hypercube} let the vector space homomorphism $\varphi$ be from Lemma \ref{L:psi}(ii). Recall that $\varphi$ is an isomorphism of vector spaces from $V_{+}$ to $\tilde{V}$.  One routinely finds that $\varphi\circ\mathbf{A}\circ\varphi^{-1}=\tilde{\mathbf{A}}$ and $\varphi\circ\mathbf{A}^{*}_{D-1}\circ\varphi^{-1}=\tilde{\mathbf{B}}$.  By this and Lemma \ref{L:ad-1} one finds the existence of an $\mathcal{A}$-module structure on $\tilde{V}$ such that the generators $x,y$ act as $\tilde{\mathbf{A}},\tilde{\mathbf{B}}$ respectively.  The $\mathcal{A}$-module structure is unique because $x,y$ generate $\mathcal{A}$.

(ii) By construction and Lemma \ref{L:psi}(ii), the vector space isomorphism $\varphi:V_{+}\to\tilde{V}$ is an isomorphism of $\mathcal{A}$-modules.  Therefore the $\mathcal{A}$-module $\varphi^{-1}(W)\subseteq V_{+}$ has diameter $\mathcal{D}-r$.  The isomorphism class of $\varphi^{-1}(W)$ and hence of $W$, can be routinely computed using the above information along with the data from Theorem \ref{T:oddcase1}. The result follows.
\hfill $\Box$ \\

%%%%%%%%%%%%%%%%%%%%%%%%%%%%%%%%%%%%%%%%%%%%%%%%%%%%%%%%%%
\section{Leonard triples from $Q_{D}$ and $\tilde{Q}_{D}$}\label{S:LTs}

In Section \ref{S:modcubes} (resp. Section \ref{S:qmodcubes}) we found that every irreducible $T$-module (resp. $\tilde{T}$-module) induces a finite-dimensional irreducible $\mathcal{A}$-module.  In this section, we show that every irreducible $T$-module (resp. $\tilde{T}$-module) induces a totally bipartite (resp. totally almost bipartite) Leonard triple of Bannai/Ito type.  We this section we describe the Leonard triples arising from our construction.

In Lemma \ref{L:sl2cube}, Go displayed an $\mathfrak{sl}_{2}$-module structure on $V$ on which $X,Y$ act as $\mathbf{A},\mathbf{A}^{*}$ respectively.  In \cite{miklavic}, Miklavi\v{c} introduced a matrix $\mathbf{A}^{\varepsilon}\in T$ that acts on $V$ as $Z$ under this $\mathfrak{sl}_{2}$-module structure.  In Section \ref{S:modcubes} we displayed an $\mathcal{A}$-module structure on $V$ (resp. $\tilde{V}$) on which $\mathbf{A},\mathbf{B}$ (resp. $\tilde{\mathbf{A}},\tilde{\mathbf{B}}$) act as $x,y$ respectively.  In this section we introduce the matrix $\mathbf{C}\in T$(resp. $\tilde{\mathbf{C}}\in\tilde{T}$) that act on $V$ (resp. $\tilde{V}$) as $z$ under these $\mathcal{A}$-module structures.

We first define the matrix $\mathbf{C}$.

\begin{definition}\label{D:IAM1}
{\normalfont With regards to Notation \ref{N:hypercube}, let $\mathbf{C}\in\mathrm{Mat}_{\mathcal{X}}(\mathbb{C})$ denote the matrix representing the action of $z$ on $V$ under the positive $\mathcal{A}$-module structure.}
\end{definition}

\begin{proposition}\label{P:zinT1}
With regards to Notation \ref{N:hypercube}, the matrix $\mathbf{C}$ from Definition \ref{D:IAM1} is contained in $T$.
\end{proposition}

\noindent {\it Proof:} 
Immediate from equation (\ref{E:rel1}) and Lemma \ref{L:ad-1}.
\hfill $\Box$ \\

\begin{definition}\label{D:WAM}
{\normalfont Given a graph $\Gamma$ with vertex set $\mathcal{X}$ and edge set $\mathcal{R}$, a matrix $\mathbf{N}\in\mathrm{Mat}_{\mathcal{X}}(\mathbb{C})$ is said to be a \emph{weighted adjacency matrix} for $\Gamma$ whenever the $(\mathbf{y},\mathbf{z})$-entry is nonzero for all adjacent $b,c\in\mathcal{X}$ and zero for all nonadjacent $\mathbf{y},\mathbf{z}\in\mathcal{X}$.}
\end{definition}

\begin{proposition}\label{P:zaction1}
With regards to Notation \ref{N:hypercube}, the matrix $\mathbf{C}$ from Definition \ref{D:IAM1} is a weighted adjacency matrix for $Q_{D}$.  For adjacent $\mathbf{y},\mathbf{z}\in\mathcal{X}$, the $(\mathbf{y},\mathbf{z})$-entry of $\mathbf{C}$ is $(-1)^{i}$ where $i$ is the minimum of $\partial(\mathbf{x},\mathbf{y})$ and $\partial(\mathbf{x},\mathbf{z})$.
\end{proposition}

\noindent {\it Proof:} 
Routine calculation using equation (\ref{E:rel1}) and Theorem \ref{T:BILT}.
\hfill $\Box$ \\

We now define the matrix $\tilde{\mathbf{C}}$.

\begin{definition}\label{D:IAM2}
{\normalfont With regards to Notation \ref{N:qhypercube}, let $\tilde{\mathbf{C}}\in\mathrm{Mat}_{\tilde{\mathcal{X}}}(\mathbb{\mathbf{C}})$ denote the matrix representing the action of $z$ on $\tilde{V}$ under the $\mathcal{A}$-module structure from Theorem \ref{T:oddcase2}.}
\end{definition}

\begin{proposition}\label{P:zinT2}
With regards to Notation \ref{N:qhypercube}, the matrix $\tilde{\mathbf{C}}$ from Definition \ref{D:IAM2} is contained in $\tilde{T}$.
\end{proposition}

\noindent {\it Proof:} 
Immediate from equation (\ref{E:rel1}) and Theorem \ref{T:oddcase2}(i).
\hfill $\Box$ \\

We have a few comments about Lemma \ref{L:Vpm}.  Recall that, when $D$ is odd, the space $V_{+}$ is a $T$-submodule of $V$.  By this and Proposition \ref{P:zinT1}, $V_{+}$ is closed under the action of $\mathbf{C}$.

\begin{proposition}\label{P:zquotient}
With regards to Notations \ref{N:hypercube}, \ref{N:qhypercube}, let $D$ be odd and let $\varphi:V\to\tilde{V}$ be from Lemma \ref{L:psi}(ii). Then $\tilde{\mathbf{C}}=\varphi\circ\mathbf{C}\circ\varphi^{-1}$ on $V$.
\end{proposition}

\noindent {\it Proof:} 
By construction, $\varphi$ is an isomorphism of $\mathcal{A}$-modules from $V_{+}$ under the positive $\mathcal{A}$-module structure to $\tilde{V}$ under the $\mathcal{A}$-module structure from Theorem \ref{T:oddcase2}.  The result follows from this information along with Definitions \ref{D:IAM1}, \ref{D:IAM2}.
\hfill $\Box$ \\

\begin{proposition}\label{P:zaction2}
With regards to Notation \ref{N:qhypercube}, the matrix $\tilde{\mathbf{C}}$ from Definition \ref{D:IAM2} is a weighted adjacency matrix for $\tilde{Q}_{D}$.  For adjacent $\mathbf{u},\mathbf{v}\in\tilde{\mathcal{X}}$, the $(\mathbf{u},\mathbf{v})$-entry of $\tilde{\mathbf{C}}$ is $(-1)^{i}$ where $i$ is the minimum of $\partial([\mathbf{x}],\mathbf{u})$ and $\partial([\mathbf{x}],\mathbf{v})$.
\end{proposition}

\noindent {\it Proof:} 
Routine calculation using Definitions \ref{D:psi}. \ref{D:IAM1}, \ref{D:IAM2} and Proposition \ref{P:zquotient}.
\hfill $\Box$ \\

From this correspondence we obtain the following results.

\begin{theorem}
With regards to Notation \ref{N:hypercube}, let $W$ denote an irreducible $T$-module.  Then the actions of $\mathbf{A},\mathbf{B},\mathbf{C}$ on $W$ form a totally bipartite Leonard triple of Bannai/Ito type.
\end{theorem}

\noindent {\it Proof:} 
Immediate from Lemma \ref{L:modLTs}.
\hfill $\Box$ \\

\begin{theorem}
With regards to Notation \ref{N:qhypercube}, let $W$ denote an irreducible $\tilde{T}$-module.  Then the actions of $\tilde{\mathbf{A}},\tilde{\mathbf{B}},\tilde{\mathbf{C}}$ on $W$ form a totally almost bipartite Leonard triple of Bannai/Ito type.
\end{theorem}

\noindent {\it Proof:} 
Immediate from Lemma \ref{L:modLTs}.
\hfill $\Box$ \\

We have found that every irreducible $T$-module (resp. $\tilde{T}$-module) induces a totally bipartite (resp. totally almost bipartite) Leonard triple of Bannai/Ito type.  Previously, we classified the totally B/AB Leonard triples of Bannai/Ito type.  In this section we describe the Leonard triples arising from our construction.

\begin{theorem}\label{T:LT1}
With regards to Notation \ref{N:hypercube}, let $D$ be even and let $\mathbf{C}$ be from Definition \ref{D:IAM1}.  Let $W$ denote an irreducible $T$-module with endpoint $r$.  Then the Leonard triple formed by the actions actions of $\mathbf{A},\mathbf{B},\mathbf{C}$ on $W$ is normalized with diameter $D-2r$.
\end{theorem}

\noindent {\it Proof:} 
By Theorem \ref{T:BILT} and Definition \ref{D:IAM1}, $W$ is an irreducible $\mathcal{A}$-module of type $B(D-2r)$ on which $x,y,z$ act as $\mathbf{A},\mathbf{B},\mathbf{C}$ respectively. By this and Proposition \ref{P:Anormalization1}, the Leonard triple formed by the actions actions of $\mathbf{A},\mathbf{B},\mathbf{C}$ on $W$ is normalized with diameter $D-2r$.
\hfill $\Box$ \\

\begin{theorem}\label{T:LT2}
With regards to Notations \ref{N:hypercube}, \ref{N:qhypercube}, let $\tilde{\mathbf{C}}$ be from Definition \ref{D:IAM2}. Let $W$ denote an irreducible $\tilde{T}$-module with endpoint $r$.  Then the Leonard triple formed by the actions of $\tilde{\mathbf{A}},\tilde{\mathbf{B}},\tilde{\mathbf{C}}$ on $W$ is $n$-normalized with diameter $\mathcal{D}-r$ where $n$ is given in the following table.
\begin{center}
\begin{tabular}{|c||c|c|}
\hline
$n$ & $\mathcal{D}$ even & $\mathcal{D}$ odd \\ \hline\hline
$r$ even & $0$ & $z$ \\ \hline
$r$ odd & $x$ & $y$ \\ \hline
\end{tabular}
\end{center}
\end{theorem}

\noindent {\it Proof:} 
By Theorem \ref{T:oddcase2} and Definition \ref{D:IAM2}, $W$ is an irreducible $\mathcal{A}$-module of type $AB(\mathcal{D}-r,n)$ on which $x,y,z$ act as $\tilde{\mathbf{A}},\tilde{\mathbf{B}},\tilde{\mathbf{C}}$ respectively. By this, Propositions \ref{P:Anormalization1} and \ref{P:nbijection}, the Leonard triple formed by the actions actions of $\mathbf{A},\mathbf{B},\mathbf{C}$ on $W$ is $n$-normalized with diameter $D-2r$.
\hfill $\Box$ \\

%%%%%%%%%%%%%%%%%%%%%%%%%%%%%%%%%%%%%%%%%%%%%%%%%%%%%%%%%%
\section{Acknowledgment}\label{S:acknowledgment}

This paper was written while the author was a graduate student at the University of Wisconsin-Madison.  The author would like to thank his advisor, Paul Terwilliger, for offering many valuable ideas and suggestions.


\small
\begin{thebibliography}{10}

%\bibitem{andrews-askey-roy}
%G. E.~Andrews, R.~Askey, and R.~Roy,
%\newblock{\em
%Special Functions,}
%\newblock{Encyclopedia of mathematics and its applications, Cambridge university press, 1999.}

\bibitem{arik-kayser}
M.~Arik, and U.~Kayserilioglu.
\newblock{The anticommutator spin algebra, its representations and quantum group invariance}
\newblock{\emph{Int. J. Mod. Phys.} A, \textbf{18} (2003), 5039-5046.}

\bibitem{bannai-ito}
E.~Bannai, and T.~Ito.
\newblock{\em Algebraic Combinatorics I: Association Schemes,}
\newblock{Benjamin \& Cummings, Mento Park, CA, 1984}.

\bibitem{bcn}
A. E.~Brouwer, A. M.~Cohen, and A.~Neumaier.
\newblock{\em Distance-Regular Graphs,}
\newblock{Springer Verlag, New York, NY, 1989}.

\bibitem{brown}
G. M. F.~Brown.
\newblock{Totally Bipartite/ABipartite Leonard triples of Bannai/Ito type}.

\bibitem{collins}
B. V. C.~Collins.
\newblock{The Terwilliger algebra of an almost-bipartite distance-regular graph and its antipodal 2-cover}
\newblock{\emph{Discrete Mathematics} \textbf{216} (2000), 35-69.}

\bibitem{curtin}
B.~Curtin.
\newblock{Modular Leonard triples}
\newblock{\emph{Linear Algebra Appl.} \textbf{424} (2007), 510-539.}

\bibitem{go}
J. T..~Go.
\newblock{The Terwilliger Algebra of the Hypercube}
\newblock{\emph{Europ. J. Combinatorics} \textbf{23} (2002), 399-429.}

%\bibitem{grillet}
%P A.~Grillet.
%\newblock{\em Abstract Algebra,}
%\newblock{Second Edition, Springer Science+Business Media, LLC, 2007}.

\bibitem{HKP}
M.~Havl\'{i}\v{c}ek, A. U.~Klimyk and S.~Po\v{s}ta,
\newblock{Representations of the cyclically symmetric $q$-deformed algebra $\mathfrak{so}_{q}(3)$,}
\newblock{\emph{Journal of Mathematical Physics}, Volume 40 Issue 4 (1990), 2135-2161.}

%\bibitem{huang}
%H.~Huang,
%\newblock{The classification of Leonard triples of QRacah type,}
%\newblock{August 3, 2011}.

\bibitem{humphries}
J. E.~Humphries.
\newblock{\em Introduction to Lie Algebras and Representation Theory.}
\newblock{Springer Science+Business Media Inc., 1972}.

%\bibitem{terw-nom}
%K.~Nomura, and P.~Terwilliger.
%\newblock{Balanced Leonard Pairs}
%\newblock{\emph{Journal of Algebra and Its Applications},}
%\newblock{June 15, 2005}.

%\bibitem{koekoek-swarttouw}
%R.~Koekoek, and R. F.~Swarttouw,
%\newblock{The Askey-scheme of hypergeometric orthogonal polynomials and its $q$-analogue, Delft University of Technoloty, Faculty of Information Technology and Systems, Department of Technical Mathematics and Informatics, report no. 98-17, 1998}
%(http://fa.its.tudelft.nl/~koekoek/askey/).


%\bibitem{passman}
%D.~Passman.
%\newblock{\em A course in ring theory.}
%\newblock{AMS Chelsea publishing, 2004}.
%%%%%%%%%%%%%%

\bibitem{ITT}
T.~Ito, K.~Tanabe, and P.~Terwilliger.
\newblock{Some algebra related to $P$- and $Q$-polynomial association schemes}
\newblock{\emph{Codes and Association Schemes} (Piscataway, NJ, 1999) Amer. Math. Soc., Providence, RI (2001), 167-192.}

\bibitem{miklavic}
\v{S}.~Miklavi\v{c}.
\newblock{Leonard triples and hypercubes,}
\newblock{\emph{Journ. Alg. Comb.}, Vol. 8 No. 3 (2008), 397-424.}

\bibitem{bigterwilliger}
P.~Terwilliger.
\newblock{The subconstituent algebra of an association scheme}
\newblock{\emph{J. Algebraic Combin.} Part I \textbf{1} (1992), 363-388; Part II \textbf{2} (1993), 73-103; Part III \textbf{2} (1993), 177-210.}

\bibitem{terwilliger}
P.~Terwilliger.
\newblock{Two linear transformations each tridiagonal With regards to an eigenbasis of the other}
\newblock{\emph{Linear Algebra Appl.} \textbf{330} (2001), 149-203.}

\bibitem{zeidler}
E.~Zeidler.
\newblock{\em Quantum Field Theory III: Gauge Theory: A Bridge between Mathematicians and Physicists.}
\newblock{Springer; 1st Edition., 2011}.








 \end{thebibliography}
\end{document}